Konrad-Zuse-Zentrum für Informationstechnik Berlin
Takustr. 7, D-14195 Berlin - Dahlem

Wolfram Koepf

Dieter Schmersau*

# Algorithms for Classical Orthogonal Polynomials

* Fachbereich Mathematik und Informatik der Freien Universität Berlin

# Algorithms for Classical Orthogonal Polynomials

Wolfram Koepf
Dieter Schmersau
koepf@zib.de

**Abstract:**

In this article explicit formulas for the recurrence equation

$$p_{n+1}(x) = (A_n\,x + B_n)\,p_n(x) - C_n\,p_{n-1}(x)$$

and the derivative rules

$$\sigma(x)\,p_n'(x) = \alpha_n\,p_{n+1}(x) + \beta_n\,p_n(x) + \gamma_n\,p_{n-1}(x)$$

and

$$\sigma(x)\,p_n'(x) = (\tilde{\alpha}_n\,x + \tilde{\beta}_n)\,p_n(x) + \tilde{\gamma}_n\,p_{n-1}(x)$$

respectively which are valid for the orthogonal polynomial solutions $p_n(x)$ of the differential equation

$$\sigma(x)\,y''(x) + \tau(x)\,y'(x) + \lambda_n\,y(x) = 0$$

of hypergeometric type are developed that depend *only* on the coefficients $\sigma(x)$ and $\tau(x)$ which themselves are polynomials w.r.t. $x$ of degrees not larger than 2 and 1, respectively.

Partial solutions of this problem had been previously published by Tricomi, and recently by Yáñez, Dehesa and Nikiforov.

Our formulas yield an algorithm with which it can be decided whether a given holonomic recurrence equation (i.e. one with polynomial coefficients) generates a family of classical orthogonal polynomials, and returns the corresponding data (density function, interval) including the standardization data in the affirmative case.

In a similar way, explicit formulas for the coefficients of the recurrence equation and the difference rule

$$\sigma(x)\,\nabla p_n(x) = \alpha_n\,p_{n+1}(x) + \beta_n\,p_n(x) + \gamma_n\,p_{n-1}(x)$$

of the classical orthogonal polynomials of a discrete variable are given that depend only on the coefficients $\sigma(x)$ and $\tau(x)$ of their difference equation

$$\sigma(x)\,\Delta\nabla y(x) + \tau(x)\,\Delta y(x) + \lambda_n\,y(x) = 0\;.$$

Here

$$\Delta y(x) = y(x+1) - y(x) \qquad \text{and} \qquad \nabla y(x) = y(x) - y(x-1)$$

denote the forward and backward difference operators, respectively. In particular this solves the corresponding inverse problem to find the classical discrete orthogonal polynomial solutions of a given holonomic recurrence equation.

## 1 Polynomials of the Hypergeometric Type

A long-standing problem in the theory of special functions whose solution can be very helpful in applied mathematics as well as in many quantum-mechanical problems of physics [18], is the determination of the differentiation formulas of the hypergeometric-type orthogonal polynomials $p_n(x)$ only from the coefficients of the differential equation

$$\sigma(x)\, y''(x) + \tau(x)\, y'(x) + \lambda_n\, y(x) = 0 \tag{1}$$

which is satisfied by these polynomials

$$y(x) = p_n(x) = k_n\, x^n + \ldots \qquad (n \in \mathbb{N}_0 := \{0, 1, 2, \ldots\}, k_n \neq 0)\ . \tag{2}$$

The coefficients $\sigma(x)$, $\tau(x)$ and $\lambda_n$ turn out to be themselves polynomials w.r.t. $x$ of degrees not larger than 2, 1 and 0, respectively.
This problem was partially solved by Tricomi ([21], Chapter IV) in the sense that he was able to calculate the coefficients $\tilde{\alpha}_n$, $\tilde{\beta}_n$ and $\tilde{\gamma}_n$ of the derivative rule

$$\sigma(x)\, p_n'(x) = (\tilde{\alpha}_n\, x + \tilde{\beta}_n)\, p_n(x) + \tilde{\gamma}_n\, p_{n-1}(x)\ . \tag{3}$$

However, his formula for $\tilde{\beta}_n$ was not only in terms of the coefficients of (1) and $k_n$, but furthermore the second highest coefficients of $p_n(x)$ were involved, and to evaluate $\tilde{\gamma}_n$, he needed to know also the coefficients of the recurrence equation

$$p_{n+1}(x) = (A_n\, x + B_n)\, p_n(x) - C_n\, p_{n-1}(x) \tag{4}$$

another structural property of orthogonal polynomial systems.
Since the polynomials $p_n(x)$ given by (2) are completely determined by the differential equation and their leading coefficients $k_n$ $(n \in \mathbb{N}_0)$, it is desirable to obtain the recurrence equation (4) *and* the derivative rule (3) from these informations alone.
Recently, Yáñez, Dehesa and Nikiforov [23] presented such formulas which, however, are additionally in terms of the constant $D_n$, given by a representation of the type

$$p_n(x) = \frac{D_n}{\rho(x)} \int\limits_C \frac{\sigma^n(s)\, \rho(s)}{(s-x)^{n+1}}\, ds \tag{5}$$

for $p_n(x)$, $\rho(x)$ being solution of the equation $(\sigma\rho)' = \tau\rho$, and $C$ being a contour satisfying certain boundary conditions. Their development is more general in the sense that they did not assume that $n$ is an integer. On the other hand, the assumption that $n$ is an integer implies that the contour $C$ is closed, the integral representation (5) being equivalent to the *Rodrigues representation*

$$p_n(x) = \frac{E_n}{\rho(x)} \frac{d^n}{dx^n} \Big(\rho(x)\, \sigma(x)^n\Big) \tag{6}$$

where

$$D_n = \frac{n!}{2\pi i}\, E_n\ , \tag{7}$$

and the solutions are classical orthogonal polynomials with density $\rho(x)$.
In this article, we represent the coefficients of both (3) and (4) in terms of $\sigma(x), \tau(x)$ and the term ratio $k_{n+1}/k_n$ alone, hence giving a complete solution of the proposed problem.

It is clear that our formulas should depend additionally on the leading coefficients $k_n$, since such a *standardization* can be prescribed arbitrarily. If one takes the *monic* standardization, i.e. $k_n \equiv 1$, then the formulas in fact depend only on the coefficients of the differential equation. For the classical orthogonal polynomials our formulas are stronger than Yáñez', Dehesa's and Nikiforov's result since $k_n$ is intrinsic part of $p_n(x)$, whereas the constants $D_n, E_n$ are not. Moreover, we will give $D_n$ and $E_n$ in terms of the coefficients of the differential equation, too. Algebraically two identities (differential equation *and* recurrence equation, e.g.) are needed to deduce the third one (derivative rule, e.g.), see [8], whereas here (kind of magic) we would like to deduce two from one. That this is possible is due to the *analytic* knowledge that orthogonal polynomial solutions of the differential equation (1) satisfy some structural properties, namely, the recurrence equation and derivative rule take special forms.
We make the general assumption that our polynomials $p_n(x)$ are orthogonal w.r.t. a measure $\mu$, i.e.

$$\int\limits_I p_n(x)\,p_m(x)\,d\mu(x) = \begin{cases} 0 & \text{if } m \neq n \\ h_n \neq 0 & \text{if } m = n \end{cases} \tag{8}$$

where $I$ denotes an appropriate integration path, for example a real interval.
Major tools in our development are the following well-known structural properties of such families of orthogonal polynomials.

**Lemma 1** Any system of polynomials $\{p_n(x) \mid n \in \mathbb{N}_0\}$, $p_n$ being of exact degree $n$, orthogonal with respect to a measure $\mu$, satisfies a three-term recurrence equation of the form (4)

$$p_{n+1}(x) = (A_n\,x + B_n)\,p_n(x) - C_n\,p_{n-1}(x) \qquad (n \in \mathbb{N}_0,\ p_{-1}(x) \equiv 0)\ ,$$

$A_n, B_n$ and $C_n$ not depending on $x$.

*Proof:* This property is well-known (see e.g. [21], Chapter IV). To prove it, one substitutes (2), equates the coefficients of $x^{n+1}$, and gets immediately that

$$A_n = \frac{k_{n+1}}{k_n}\ . \tag{9}$$

With this choice, we study the difference $p_{n+1}(x) - A_n\,x\,p_n(x)$. Since this is a polynomial of degree not larger than $n$, it can be decomposed as

$$p_{n+1}(x) - A_n\,x\,p_n(x) = \sum_{j=0}^{n} d_j\,p_n(x)\ .$$

We choose $m \leq n-2$ and multiply by $p_m(x)$. Integrating with respect to $\mu$ yields

$$\int\limits_I p_m(x)\,p_{n+1}(x)\,d\mu(x) - \int\limits_I p_m(x)\,A_n\,x\,p_n(x)\,d\mu(x) = d_m\,h_m$$

where on the right hand side (8) was applied. Both left hand integrals vanish since $p_m(x)$ is orthogonal to $p_{n+1}(x)$ and since $xp_m(x)$ as a polynomial of degree not larger than $n-1$ is orthogonal to $p_n(x)$, implying $d_m = 0$. This gives the result. □

The second important structural property for our considerations is given by **Lemma 2** Any system of polynomials $\{p_n(x) \mid n \in \mathbb{N}_0\}$, $p_n$ being of exact degree $n$, that are solutions of the differential equations (1) and furthermore orthogonal with respect to a measure $\mu(x) = \rho(x)\,dx$ having weight function $\rho(x) \geq 0$, satisfies a derivative rule of the form (3)

$$\sigma(x)\,p_n'(x) = (\tilde{\alpha}_n\,x + \tilde{\beta}_n)\,p_n(x) + \tilde{\gamma}_n\,p_{n-1}(x) \qquad (n \in \mathbb{N} := \{1,2,3,\ldots\})\ ,$$

$\tilde{\alpha}_n, \tilde{\beta}_n$ and $\tilde{\gamma}_n$ not depending on $x$.

*Proof:* Substituting (2) and equating the coefficients of $x^{n+1}$, one gets immediately that

$$\tilde{\alpha}_n = a\,n\ . \tag{10}$$

In [13], § 5 it is shown by an elementary argument that under the given conditions the solutions $p_n(x)$ of the differential equations (1) are orthogonal with respect to the weight function

$$\rho(x) := \frac{C}{\sigma(x)}\,e^{\int \frac{\tau(x)}{\sigma(x)}\,dx} \geq 0\ , \tag{11}$$

given by *Pearson's differential equation*

$$\frac{d}{dx}\Big(\sigma(x)\,\rho(x)\Big) = \tau(x)\,\rho(x)$$

for a suitable constant $C$, in a suitable interval $I$ (depending on the zeros of $\sigma(x)$). Hence multiplying (1) by $\rho(x)$, the differential equation takes the selfadjoint form

$$\frac{d}{dx}\Big(\sigma(x)\,\rho(x)\,y'(x)\Big) + \lambda_n\,\rho(x)\,y(x) = 0\ .$$

Using this identity, Tricomi showed that ([21], IV (4.10))

$$\int_I \sigma(x)\,\rho(x)\,p_n'(x)\,f(x)\,dx = 0 \tag{12}$$

for any polynomial $f(x)$ of degree $\leq n-2$. If (10) holds, then the degree of $\sigma(x)\,p_n'(x) - \tilde{\alpha}_n\,x$ is $\leq n$. Hence one can write

$$\sigma(x)\,p_n'(x) - \tilde{\alpha}_n\,x = \sum_{j=0}^{n} e_j\,p_n(x)\ .$$

As above, from (12) one can deduce that $e_j = 0$ for $0 \leq j \leq n-2$ (see [21], Chapter IV). $\square$

An immediate consequence is the following

**Corollary 1** Any system of polynomials $\{p_n(x) \mid n \in \mathbb{N}_0\}$, $p_n$ being of exact degree $n$, that are solutions of the differential equation (1) and furthermore orthogonal with respect to a measure $\mu(x) = \rho(x)\,dx$ having weight function $\rho(x) \geq 0$, satisfies a derivative rule of the form

$$\sigma(x)\,p_n'(x) = \alpha_n\,p_{n+1}(x) + \beta_n\,p_n(x) + \gamma_n\,p_{n-1}(x) \qquad (n \in \mathbb{N})\ , \tag{13}$$

$\alpha_n, \beta_n$ and $\gamma_n$ not depending on $x$.

*Proof:* Substituting (2) in (13), and equating the coefficients of $x^{n+1}$, one gets immediately that

$$\alpha_n = a\,n\,\frac{k_n}{k_{n+1}}\ . \tag{14}$$

Substituting (4) in (3) one gets moreover

$$\begin{aligned}
\sigma(x)\,p_n'(x) &= (\tilde{\alpha}_n\,x + \tilde{\beta}_n)\,p_n(x) + \tilde{\gamma}_n\,p_{n-1}(x) \\
&= \left(\frac{\tilde{\alpha}_n}{A_n}\left(p_{n+1}(x) - B_n\,p_n(x) + C_n\,p_{n-1}(x)\right) + \tilde{\beta}_n\right)\,p_n(x) + \tilde{\gamma}_n\,p_{n-1}(x)\ ,
\end{aligned}$$

hence (13) is valid with

$$\alpha_n = \frac{\tilde{\alpha}_n}{A_n}\ , \qquad \beta_n = \tilde{\beta}_n - \tilde{\alpha}_n\frac{B_n}{A_n}\ , \qquad \gamma_n = \tilde{\gamma}_n + \tilde{\alpha}_n\frac{C_n}{A_n}\ . \qquad \square$$

## 2 Classical Orthogonal Polynomials of an Interval

In this section we give the proposed explicit recurrence equation and derivative rule formulas. Assume a family of differential equations (1) is given for $n \in \mathbb{N}_0$, with continuous functions $\sigma(x), \tau(x)$, and constants $\lambda_n$, and we search for polynomial solutions (2) of degree $n$. Then since $p_1(x)$ is linear, one deduces that $\tau(x)$ must be an at most linear polynomial, and since $p_2(x)$ is quadratic, one deduces that $\sigma(x)$ must be an at most quadratic polynomial [3]. Hence we may assume that

$$\sigma(x) := ax^2 + bx + c\ , \qquad \tau(x) := dx + e\ . \tag{15}$$

Equating coefficients of the highest powers $x^n$ in (1) for generic $p_n(x)$, given by (2), one deduces that moreover

$$an(n-1) + dn + \lambda_n = 0 \qquad \text{or} \qquad \lambda_n = -(an(n-1) + dn)\ . \tag{16}$$

Hence only if the differential equation takes the special form

$$(ax^2 + bx + c)\,y''(x) + (dx + e)\,y'(x) - (an(n-1) + dn)\,y(x) = 0\ , \tag{17}$$

it can have polynomial solutions.
Moreover we can assume that $\lambda_n \neq 0$ for $n \in \mathbb{N}$, hence $a(n-1) + d \neq 0$ for $n \in \mathbb{N}$ since otherwise no orthogonal polynomial solutions can exist. This is discussed in detail in [13]. In particular, $d \neq 0$.
In the following theorem, we give explicit representations of the corresponding recurrence equation and derivative rule in terms of the given $a, b, c, d, e$ and the term ratio $k_{n+1}/k_n$.

**Theorem 1** Let $p_n(x) = k_n\,x^n + \ldots$ $(n \in \mathbb{N}_0)$ be a family of polynomial solutions of the system of differential equations (17) that are orthogonal with respect to a weight function $\rho(x)$. Then the derivative rule (13)

$$\sigma(x)\,p_n'(x) = \alpha_n\,p_{n+1}(x) + \beta_n\,p_n(x) + \gamma_n\,p_{n-1}(x)$$

is valid with

$$\alpha_n = a\,n\,\frac{k_n}{k_{n+1}}$$

$$\beta_n = \frac{n\,(a\,(n-1)+d)\,(b\,d-2\,a\,e)}{(2\,a\,(n-1)+d)\,(2\,a\,n+d)} \tag{18}$$

$$\gamma_n = \frac{n\,(a(n-1)+d)(a(n-2)+d)(n\,(an+d)\,(4\,a\,c-b^2)+a\,e^2+c\,d^2-b\,d\,e)}{(a\,(2n-1)+d)\,(a(2n-3)+d)\,(2a(n-1)+d)^2}\,\frac{k_n}{k_{n-1}}\ , \tag{19}$$

and the recurrence equation (4)

$$p_{n+1}(x) = (A_n\,x + B_n)\,p_n(x) - C_n\,p_{n-1}(x)$$

is valid with

$$A_n = \frac{k_{n+1}}{k_n}$$

$$B_n = \frac{k_{n+1}}{k_n}\cdot\frac{2bn\,(a(n-1)+d)+e(d-2a)}{(2a(n-1)+d)\,(2an+d)} \tag{20}$$

and

$$C_n = -\frac{k_{n+1}}{k_n}\,\frac{\gamma_n}{a(n-1)+d}\ , \tag{21}$$

$\gamma_n$ being given by (19).

*Proof:* The values of $A_n$ and $\alpha_n$ were already obtained in Lemma 1 and Corollary 1. By Lemma 1 the polynomials satisfy a recurrence equation of type (4):

$$p_{n+1}(x) = (A_n\,x + B_n)\,p_n(x) - C_n\,p_{n-1}(x)\ . \tag{22}$$

Next, we differentiate (22) twice and get

$$p'_{n+1}(x) = A_n\,p_n(x) + (A_n\,x + B_n)\,p'_n(x) - C_n\,p'_{n-1}(x) \tag{23}$$

and

$$p''_{n+1}(x) = 2A_n\,p'_n(x) + (A_n\,x + B_n)\,p''_n(x) - C_n\,p''_{n-1}(x)\ .$$

We multiply the last equation by $\sigma(x)$

$$\sigma(x)\,p''_{n+1}(x) = 2A_n\,\sigma(x)\,p'_n(x) + (A_n\,x + B_n)\,\sigma(x)\,p''_n(x) - C_n\,\sigma(x)\,p''_{n-1}(x)$$

and use the differential equation to replace the second derivatives by those of lower order

$$\begin{aligned} -\Big(\tau(x)\,p'_{n+1}(x) + \lambda_{n+1}\,p_{n+1}(x)\Big) \;&=\; 2A_n\,\sigma(x)\,p'_n(x) - (A_n\,x + B_n)\,\Big(\tau(x)\,p'_n(x) + \lambda_n\,p_n(x)\Big) \\ &\quad + C_n\,\Big(\tau(x)\,p'_{n-1}(x) + \lambda_{n-1}\,p_{n-1}(x)\Big)\ . \end{aligned}$$

After substituting (23) on the left hand side, and subtracting $\tau\,(Ax+B)\,p'_n - \tau\,C_n\,p'_{n-1}$, we arrive at

$$-\tau(x)\,A_n\,p_n(x) - \lambda_{n+1}\,p_{n+1}(x) = 2A_n\,\sigma(x)\,p'_n(x) - (A_n\,x + B_n)\,\lambda_n\,p_n(x) + C_n\,\lambda_{n-1}\,p_{n-1}(x)\ .$$

Next, on the right hand side, we replace $(A_n\,x + B_n)\,p_n(x)$ by $p_{n+1}(x) + C_n\,p_{n-1}(x)$ according to (22), and get

$$-\tau(x)\,A_n\,p_n(x) - \lambda_{n+1}\,p_{n+1}(x) = 2A_n\,\sigma(x)\,p'_n(x) - \lambda_n\,p_{n+1}(x) + (\lambda_{n-1} - \lambda_n)\,C_n\,p_{n-1}(x)\ ,$$

or rewritten

$$(\lambda_n - \lambda_{n+1})\, p_{n+1}(x) = A_n\, (2\,\sigma(x)\, p_n'(x) + \tau(x)\, p_n(x)) + (\lambda_{n-1} - \lambda_n)\, C_n\, p_{n-1}(x)\ .$$

Now we substitute the representation of Corollary 1

$$\sigma(x)\, p_n'(x) = \alpha_n\, p_{n+1}(x) + \beta_n\, p_n(x) + \gamma_n\, p_{n-1}(x)$$

to deduce

$$\Big((\lambda_n - \lambda_{n+1}) - 2\, A_n\, \alpha_n\Big)\, p_{n+1}(x) = A_n\, (2\,\beta_n + \tau(x))\, p_n(x) + \Big(2\, A_n\, \gamma_n + (\lambda_{n-1} - \lambda_n)\, C_n\Big)\, p_{n-1}(x)$$

after subtracting $2\, A_n\, \alpha_n\, p_{n+1}(x)$. Replacing $p_{n+1}(x)$ according to (22), we arrive at the identity

$$\Big(A_n\, (2\,\beta_n + \tau(x)) + (2\, A_n\, \alpha_n - (\lambda_n - \lambda_{n+1}))(A_n\, x + B_n)\Big)\, p_n(x) =$$
$$-\Big(((\lambda_n - \lambda_{n+1}) - 2\, A_n\, \alpha_n)\, C_n + 2\, A_n\, \gamma_n + (\lambda_{n-1} - \lambda_n)\, C_n\Big)\, p_{n-1}(x)\ .$$

Since $p_n(x)$ is a polynomial of exact degree $n$, this relation can only be valid if the coefficients of both $p_n(x)$ and $p_{n-1}(x)$ vanish, since otherwise the polynomial on the left hand side has degree $\geq n$ whereas the polynomial on the right hand side has degree $n-1$, a contradiction. The coefficients must vanish as polynomials in $x$, and equating coefficients we are led to the three equations

$$A_n\, (\lambda_{n+1} - \lambda_n + d + 2\,\alpha_n\, A_n) = 0\ .$$

$$A_n\, e + 2\, A_n\, B_n\, \alpha_n + 2\, A_n\, \beta_n - B_n\, \lambda_n + B_n\, \lambda_{n+1} = 0$$

and

$$2\, A_n\, C_n\, \alpha_n - 2\, A_n\, \gamma_n - C_n\, \lambda_{n-1} + C_n\, \lambda_{n+1} = 0\ .$$

Whereas the first of these equation does not contain any news but restates a relationship between $A_n$, $\alpha_n$ and $\lambda_n$, the second and third of these equations (using (9), (14) and (16)) can be rewritten as

$$B_n = \frac{(e + 2\,\beta_n)}{d}\, \frac{k_{n+1}}{k_n} \tag{24}$$

and as (21). Hence $B_n$ and $C_n$ are known as soon as $\beta_n$ and $\gamma_n$ are.
We finally need two more equations to find $\beta_n$ and $\gamma_n$. To deduce one of these equations, and to find $\beta_n$, and hence $B_n$, we substitute

$$p_n(x) = k_n\, x^n + k_n'\, x^{n-1} + k_n''\, x^{n-2} + \ldots \tag{25}$$

in the three equations considered, namely the differential equation, the recurrence equation and the derivative rule. As we already saw, equating the coefficients of the highest powers of $x$ yields (16), (9) and (14). If we equate the coefficients of the next highest powers of $x$, we get three more equations, involving two more variables though, namely $k_n'$ and $k_{n+1}'$. These are the equations

$$b\, n\, k_n - e\, n\, k_n - b\, n^2\, k_n - 2\, a\, k_n' + d\, k_n' + 2\, a\, n\, k_n' = 0\ , \tag{26}$$

$$B_n\, k_n + A_n\, k_n' - k_{n+1}' = 0\ , \tag{27}$$

and

$$b\, n\, k_n - \beta_n\, k_n - a\, k_n' + a\, n\, k_n' - \alpha_n\, k_{n+1}' = 0\ . \tag{28}$$

Equation (26) immediately gives

$$\frac{k_n'}{k_n} = \frac{n\ (b\,(n-1)+e)}{2\,a\,(n-1)+d}\ , \tag{29}$$

whereas from (27)–(28) one can eliminate $k_{n+1}'$. This gives a second equation between $B_n$ and $\beta_n$ which together with (24) and (29) yields (18) and (20).
To deduce $\gamma_n$, we equate the coefficients of the next highest powers in the differential equation, recurrence equation and derivative rule, introducing two more auxiliary variables $k_n''$ and $k_{n+1}''$ which can be eliminated. This procedure generates one more equation between $C_n$ and $\gamma_n$ finally deducing (19). □

Note that the results given in Theorem 1 can also be deduced completely automatically by elimination methods based on Gröbner basis calculations. With the computer algebra systems Maple and REDUCE we were successful doing so. For the purpose of finding $A_n, B_n, C_n$, we substitute (25) in the differential equations for $p_n(x)$ and for $p_{n+1}(x)$, and in the recurrence equation. Equating the three highest coefficients in any of these three equations yields nine nonlinear equations in the nine unknowns

$$A_n, B_n, C_n, \lambda_n, \lambda_{n+1}, k_n', k_{n+1}', k_n'', k_{n+1}''\ .$$

By a Gröbner basis computation (invoked by the `solve` command of the utilized computer algebra system) it turns out that there is a unique solution, given by Theorem 1, see also Corollary 4 and (33). Note that therefore the formulas for $A_n, B_n$, and $C_n$ of Theorem 1 are valid without the hypothesis of a weight function $\rho(x)$.
Similarly, to find $\alpha_n, \beta_n, \gamma_n$, we substitute (25) in the differential equations for $p_n(x)$ and for $p_{n+1}(x)$, and in the derivative rule. Equating the three highest coefficients in any of these three equations yields nine nonlinear equations in the nine unknowns

$$\alpha_n, \beta_n, \gamma_n, \lambda_n, \lambda_{n+1}, k_n', k_{n+1}', k_n'', k_{n+1}''\ ,$$

and a Gröbner basis computation generates the unique solution, given by Theorem 1. Note that we were not able to separate the two problems in a similar way based on hand calculations.

Our theorem has immediate consequences.

**Corollary 2** Let $p_n(x) = k_n\, x^n + \ldots\ (n \in \mathbb{N}_0)$ be a family of polynomial solutions of the system of differential equations (17) that are orthogonal with respect to a weight function $\rho(x)$. Then the derivative rule (3)

$$\sigma(x)\, p_n'(x) = (\tilde{\alpha}_n\, x + \tilde{\beta}_n)\, p_n(x) + \tilde{\gamma}_n\, p_{n-1}(x)$$

is valid with

$$\tilde{\alpha}_n = a\, n$$

$$\tilde{\beta}_n = \frac{(a\,b\,(n-1) - a\,e + b\,d)\,n}{2\,a(n-1)+d} \tag{30}$$

$$\tilde{\gamma}_n = \gamma_n - \alpha_n\, C_n = \frac{a\,(2n-1)+d}{a\,(n-1)+d}\,\gamma_n\ , \tag{31}$$

$\alpha_n$, $\gamma_n$ and $C_n$ being given by Theorem 1.

*Proof:* Substituting (4) in (13) yields (3) with

$$\tilde{\alpha}_n = \alpha_n\,A_n\;, \qquad \tilde{\beta}_n = \alpha_n\,B_n + \beta_n\;, \qquad \tilde{\gamma}_n = \gamma_n - \alpha_n\,C_n\;.$$

This yields the result. □

Note that Theorem 1 describes the variety of different recurrence equation formulas known in the literature ([1], 22.7) by *one single formula*. Similarly all the different derivative rule formulas ([1], 22.8) are governed by a single formula through Corollary 2.
Theorem 1 shows in particular that the recurrence equation and derivative rule can be obtained by purely rational arithmetic whenever

$$\frac{k_{n+1}}{k_n} \in \mathbb{Q}(n)\;,$$

i.e., if $k_n$ is a *hypergeometric term*. This is obviously true if $k_n \equiv 1$, i.e., in the monic case. But also all other standardizations that are used in practice (see e.g. [1], Chapter 22) are of this type.[1]
In the case of the *orthonormal* standardization given by

$$h_n \equiv 1$$

it is not in general true that $k_n$ is a hypergeometric term. On the other hand, if $k_n$ is a hypergeometric term, $h_n$ inherits this property.

**Corollary 3** Let $p_n(x) = k_n\,x^n + \ldots$ $(n \in \mathbb{N}_0)$ be a family of orthogonal polynomial solutions of the system of differential equations (17). Then the relation

$$\frac{h_{n+1}}{h_n} = \frac{(n+1)\,(an+d)(a(n-1)+d)}{(a(2n+3)+d)(a(2n+1)+d)}\left(c+\frac{b(n+1)+e}{(2an+d)^2}\Big((ae-bd)-abn\Big)\right)\cdot\left(\frac{k_{n+1}}{k_n}\right)^2 \qquad (32)$$

is valid.

*Proof:* Tricomi ([21], IV (2.2), see also [1], (22.1.5)) proved that

$$C_n = \frac{A_n}{A_{n-1}}\,\frac{h_n}{h_{n-1}}\;.$$

An application of Theorem 1 yields (32). □

Next we would like to give a general formula for the term ratio of the coefficients $k_n'$ in terms of the given term ratio of $k_n$.

**Corollary 4** Let $p_n(x) = k_n\,x^n + k_n'\,x^{n-1} + \ldots$ $(n \in \mathbb{N}_0)$ be a family of orthogonal polynomial solutions of the system of differential equations (17). Then the relation

$$\frac{k_{n+1}'}{k_n'} = \frac{n+1}{n}\cdot\frac{(bn+e)(2a(n-1)+d)}{(b(n-1)+e)(2an+d)}\cdot\frac{k_{n+1}}{k_n}$$

is valid.

[1] Only in one instance, this is not so: For the Chebyshev polynomials $T_n(x)$ one has $k_{n+1}/k_n = 2$ $(n \in \mathbb{N})$, and $k_1/k_0 = 1$. If one redefines $T_0(x) := 1/2$, then $k_{n+1}/k_n \equiv 2 \in \mathbb{Q}(n)$.

*Proof:* This follows from

$$\frac{k'_{n+1}}{k'_n} = \frac{k'_{n+1}}{k_{n+1}} \cdot \frac{k_{n+1}}{k_n} \cdot \frac{k_n}{k'_n}$$

using (29). □

Similarly, one gets

$$\frac{k''_n}{k_n} = \frac{n\,(n-1)\,(n^2\,b^2 - 3\,n\,b^2 + 2\,n\,b\,e + 2\,c\,n\,a - 2\,c\,a - 3\,b\,e + 2\,b^2 + e^2 + c\,d)}{2\,(2\,a\,n - 2\,a + d)\,(d - 3\,a + 2\,a\,n)}\,, \tag{33}$$

also deduced by the automatic elimination method mentioned before, and a similar equation for $k''_{n+1}/k''_n$, see the Appendix.

We furthermore obtain the term ratio of the numbers $D_n$ of an integral representation of type (5) considered in [23]. Note that in the orthogonal polynomial case the contour $C$ is closed, and hence by Cauchy's integral formula representation (5) is equivalent to a *Rodrigues representation* (6) with

$$D_n = \frac{n!}{2\pi i}\,E_n\,.$$

We get for $E_n$ and $D_n$, respectively

**Corollary 5** Let $p_n(x) = k_n\,x^n + \ldots$ $(n \in \mathbb{N}_0)$ be a family of orthogonal polynomial solutions of the system of differential equations (17). Then $p_n(x)$ have a Rodrigues representation (6) and an integral representation (5) with closed contour $C$ surrounding $s = x$, and one has for $E_n$ and $D_n$ the relations

$$\frac{D_{n+1}}{D_n} = \frac{(n+1)\,(a(n-1)+d)}{(a(2n-1)+d)(2an+d)} \cdot \frac{k_{n+1}}{k_n}\,,$$

and

$$\frac{E_{n+1}}{E_n} = \frac{(a(n-1)+d)}{(a(2n-1)+d)(2an+d)} \cdot \frac{k_{n+1}}{k_n}\,.$$

*Proof:* In ([23], (13)) it was shown that

$$\frac{1}{A_n} = \frac{D_n}{D_{n+1}} \cdot \frac{(n+1)\,(a(n-1)+d)}{(a(2n-1)+d)\,(2an+d)}\,.$$

An application of (9) leads to the term ratio for $D_n$. The term ratio for $E_n$ follows then from (7). □

Note that Corollary 5 again describes all the different Rodrigues formulas ([1], 22.11) known in the literature by one single formula.

It is well-known ([3], see also [4], [13]) that polynomial solutions of (1) can be classified according to the zeros of $\sigma(x)$, leading to the normal forms of Table 1 besides linear transformations $x \mapsto Ax + B$. The type of differential equation that we consider is invariant under such a transformation. Orthogonal polynomial solutions according to this classification exist if and only if the function

$$\rho(x) = \frac{C}{\sigma(x)}\,e^{\int \frac{\tau(x)}{\sigma(x)}\,dx}$$

given by (11) yields a weight function in the interval given by the zeros of $\sigma(x)$, i.e. the corresponding integrals converge and $\rho(x) \geq 0$ for some $C$.

1. $a=b=c=e=0,\ d=1 \quad\Longrightarrow\quad p_n(x)=x^n$,
2. $a=b=e=0,\ c=1,\ d=-2 \quad\Longrightarrow\quad p_n(x)=H_n(x)$, the *Hermite polynomials*,
3. $a{=}c{=}0,\ b{=}1,\ d{=}-1,\ e{=}\alpha{+}1 \quad\Longrightarrow\quad p_n(x)=L_n^{(\alpha)}(x)$, the *Laguerre polynomials*,

4a. $a=1,\ b=c=d=e=0, \quad\Longrightarrow\quad p_n(x)=x^n$,

4b. $a=1,\ b=c=0,\ d=\alpha+2,\ e=2 \quad\Longrightarrow\quad p_n(x)=B_n^{(\alpha)}(x)$, the *Bessel polynomials*,

5. $a{=}1, b{=}0, c{=}-1, d{=}\alpha{+}\beta{+}2, e{=}\alpha{-}\beta \;\Rightarrow\; p_n(x){=}P_n^{(\alpha,\beta)}(x)$, the *Jacobi polynomials*.

Table 1: Normal Forms of Polynomial Solutions

This shows that the only orthogonal polynomial solutions are linear transforms of the Hermite, Laguerre, and Jacobi polynomials, hence using a mathematical dictionary one can always deduce the recurrence equation and derivative rules. Note, however, that this approach (in general) requires the work with radicals, namely the zeros of the quadratic polynomial $\sigma(x)$, whereas our approach is completely rational: Given $k_{n+1}/k_n \in \mathbb{Q}(n)$, the recurrence equation and derivative rules are given rationally in Theorem 1.
Note that the formulas of Theorem 1 are also valid for the Bessel polynomials ([18], p. 24)

$$B_n^{(\alpha)}(x) = \frac{(2n)!\,x^n}{n!\,2^n}\,{}_1F_1\left(\left.\begin{array}{c}-n\\ -2n\end{array}\right|\frac{2}{x}\right) = \frac{e^{2/x}}{2^n}\frac{d^n}{dx^n}\left(x^{2n}\,e^{-2/x}\right)\ .$$

This is so since the Bessel polynomials *do* satisfy both a recurrence equation and a derivative rule of the desired type (see e.g. [23]), despite the fact that the corresponding function

$$\rho(x) = \frac{C}{\sigma(x)}\,e^{\int\frac{\tau(x)}{\sigma(x)}\,dx} = C\,x^\alpha\,e^{-2/x}$$

does not constitute a weight function on the real axis. The validity of both a recurrence equation and a derivative rule of the given types, however, was the only assumption in the proof of Theorem 1.
Although the Jacobi polynomials $P_n^{(\alpha,\beta)}(x)$ do only constitute orthogonal polynomials for $\alpha,\beta>-1$, by a simple argument it can be shown that the structural properties like recurrence equation and derivative rule remain valid for arbitrary values of $\alpha,\beta$. A similar comment applies to the other parameterized families of Table 1. Hence Theorem 1 is valid also in these cases.
Theorem 1 is even valid in the case of Table 1:4a, and its recurrence equation part also for Table 1:1 with the trivial solution $p_n(x)=x^n$. In both cases we have the recurrence equation $p_{n+1}(x)=x\,p_n(x)$, and in the first case we receive the derivative rule $x^2\,p_n'(x)=n\,p_{n+1}(x)$. Note that there is another derivative rule $x\,p_n'(x)=n\,p_n(x)$ which cannot be discovered by Theorem 1.

In the next section we will use the fact that these equations are given explicitly to solve an inverse problem.

## 3 The Inverse Characterization Problem

Assume you have a polynomial system given by a differential equation (1). Then by the classification of Table 1 it is easy to identify the system. On the other hand, given an arbitrary holonomic three-term recurrence equation

$$q_n(x)\,P_{n+2}(x) + r_n(x)\,P_{n+1}(x) + s_n(x)\,P_n(x) = 0 \qquad (q_n(x), r_n(x), s_n(x) \in \mathbb{Q}[n,x])\ , \tag{34}$$

it is less obvious to find out whether there is a polynomial system

$$P_n(x) = k_n\,x^n + \ldots \qquad (n \in \mathbb{N}_0\,, k_n \neq 0)$$

satisfying (34), being a linear transform of one of the classical systems (Hermite, Laguerre, Jacobi, Bessel), and to identify the system in the affirmative case. In this section we present an algorithm for this purpose. Note that Koornwinder and Swarttouw have also considered this question and propose a solution based on the careful ad hoc analysis of the input polynomials $q_n, r_n$, and $s_n$. Their Maple implementation works for a part of the so-called Askey-Wilson scheme ([2], see also [12]).

Let us start with a recurrence equation of type (34). We assume that neither $q_{n-1}(x)$ nor $s_n(x)$ has a nonnegative integer zero since otherwise this recurrence equation cannot be used to determine $P_n(x)$ iteratively from $P_0(x)$ (with $P_{-1}(x) \equiv 0$) for all $n \geq 1$ or is worthless in the backward direction. Define

$$N := \begin{cases} 0, & \text{if neither } q_{n-1}(x) \text{ nor } s_n(x) \text{ have a nonnegative integer zero} \\ \max\{n \in \mathbb{N}_0 \mid n \text{ is a zero of either } q_{n-1}(x) \text{ or } s_n(x)\} + 1, & \text{otherwise} \end{cases} .$$

Then we consider $p_n(x) := P_{n+N}(x)$ instead of $P_n(x)$, see § 4 for an example with $N > 0$. In this situation we rewrite (34) by substituting $n$ by $n+N$ and replacing $P_n(x)$ by $p_{n-N}(x)$. For simplicity we rename $q_n(x), r_n(x)$ and $s_n(x)$, and assume in the sequel that the recurrence equation

$$q_n(x)\,p_{n+2}(x) + r_n(x)\,p_{n+1}(x) + s_n(x)\,p_n(x) = 0 \qquad (q_n(x), r_n(x), s_n(x) \in \mathbb{Q}[n,x]) \tag{35}$$

is valid, but now neither $q_{n-1}(x)$ nor $s_n(x)$ have nonnegative integer zeros. We search for solutions

$$p_n(x) = k_n\,x^n + \ldots \qquad (n \in \mathbb{N}_0\,, k_n \neq 0) \tag{36}$$

which reads in terms of the original family $P_n$ $(n = N, N+1, \ldots)$

$$P_{n+N}(x) = k_n\,x^n + \ldots \qquad (n \in \mathbb{N}_0)\ .$$

Next, we divide (35) by $q_n(x)$, and replace $n$ by $n-1$. This brings (35) into the form

$$p_{n+1}(x) = t_n(x)\,p_n(x) + u_n(x)\,p_{n-1}(x) \qquad (t_n(x), u_n(x) \in \mathbb{Q}(n,x))\ . \tag{37}$$

For $p_n(x)$ being a linear transform of a classical orthogonal system, there is a recurrence equation (4)

$$p_{n+1}(x) = (A_n\,x + B_n)\,p_n(x) - C_n\,p_{n-1}(x) \qquad (A_n, B_n, C_n \in \mathbb{Q}(n),\ A_n \neq 0)\ , \tag{38}$$

therefore (37) and (38) must agree. We would like to conclude that $t_n(x) = A_n\, x + B_n$, and $u_n(x) = -C_n$ which follows if we can show that $p_n(x)/p_{n-1}(x) \notin \mathbb{Q}(n,x)$. To prove this assertion, we assume that

$$\frac{p_n(x)}{p_{n-1}(x)} \in \mathbb{Q}(n,x) \ .$$

Hence there are $P(n,x) \in \mathbb{Q}[n][x]$ and $Q(n,x) \in \mathbb{Q}[n][x]$ with $\gcd_x(P(n,x), Q(n,x)) = 1$ such that the relation

$$Q(n,x)\, p_n(x) = P(n,x)\, p_{n-1}(x) \tag{39}$$

holds. It is a classical result for orthogonal polynomials that $\gcd_x(p_n(x), p_{n-1}(x)) = 1$ since their zeros separate each other (see e.g. [21], IV.6). Hence from (39) we conclude that

$$P(n,x) = S_n\, p_n(x) \qquad \text{and} \qquad Q(n,x) = S_n\, p_{n-1}(x) \ .$$

Since by assumption $P(n,x) \in \mathbb{Q}[n][x]$ should be a polynomial of fixed degree with respect to $x$, and since $p_n(x)$ has degree $n$, this gives an obvious contradiction.
Therefore we can conclude that $t_n(x) = A_n\, x + B_n$, and $u_n(x) = -C_n$. Hence if (37) does not have this form, i.e., if either $t_n(x)$ is not linear in $x$ or $u_n(x)$ is not a constant with respect to $x$, we see that $p_n(x)$ cannot be a linear transform of a classical orthogonal polynomial system. In the positive case, we can assume the form (38).
Since we propose solutions (36), equating the coefficients of $x^{n+1}$ in (38) we get

$$\frac{k_{n+1}}{k_n} = A_n = \frac{v_n}{w_n} \qquad (v_n, w_n \in \mathbb{Q}[n]) \ . \tag{40}$$

Hence the given $A_n = v_n/w_n \in \mathbb{Q}(n)$ generates the term ratio $k_{n+1}/k_n$, and in particular $k_n$ turns out to be a hypergeometric term which is uniquely determined by (40) up to a normalization constant $k_0 = p_0(x)$. Since the zeros of $w_n$ correspond to the zeros of $q_{n-1}(x)$, $k_n$ is defined by (40) for all $n \in \mathbb{N}$ from $k_0$.
In the next step we can eliminate the dependency of $k_n$ by generating a recurrence equation for the corresponding monic polynomials $\tilde{p}_n(x) = p_n(x)/k_n$. For $\tilde{p}_n(x)$ we get by (40)

$$\tilde{p}_{n+1}(x) = \left(x + \frac{B_n}{A_n}\right) \tilde{p}_n(x) - \frac{C_n}{A_n\, A_{n-1}}\, \tilde{p}_{n-1}(x) = \left(x + \tilde{B}_n\right) \tilde{p}_n(x) - \tilde{C}_n\, \tilde{p}_{n-1}(x)$$

with

$$\tilde{B}_n = \frac{B_n}{A_n} \in \mathbb{Q}(n) \qquad \text{and} \qquad \tilde{C}_n = \frac{C_n}{A_n\, A_{n-1}} \in \mathbb{Q}(n) \ .$$

Then our formulas (20)–(21) read in terms of $\tilde{B}_n$ and $\tilde{C}_n$

$$\tilde{B}_n = \frac{2bn\,(a(n-1)+d) + e(d-2a)}{(2a(n-1)+d)\,(2an+d)} \tag{41}$$

and

$$\tilde{C}_n = \frac{-n\,(a(n-2)+d)}{(a\,(2n-1)+d)(a\,(2n-3)+d)} \left(c + \frac{b(n-1)+e}{(2a\,(n-1)+d)^2}\Big((ae-bd) - ab\,(n-1)\Big)\right), \tag{42}$$

and these are independent of $k_n$ by construction.
Now we would like to deduce $a, b, c, d$ and $e$ from (41)–(42). Note that as soon as we have found these five values, we can apply a linear transform (according to the zeros of $\sigma(x)$) to bring the differential equation in one of the forms of Table 1 which finally gives us the desired information.
We can assume that $\tilde{B}_n$ and $\tilde{C}_n$ are in lowest terms. If the degree of either the numerator or the denominator of $\tilde{B}_n$ is larger than 2, then by (41) $p_n(x)$ is not a classical system. Similarly, if the degree of either the numerator or the denominator of $\tilde{C}_n$ is larger than 4, by (42) the same conclusion follows.
Otherwise we can multiply (41) and (42) by their common denominators, and bring them therefore in polynomial form. Both resulting equations must be polynomial identitites in the variable $n$, hence all of their coefficients must vanish. This gives a nonlinear system of equations for the unknowns $a, b, c, d$ and $e$. Any solution of this system with not both $a$ and $d$ being zero yields a differential equation (17), and hence given such a solution one can decide whether the corresponding solutions $p_n(x)$ are generated by a density (11). Therefore our question can be resolved in this case.
If the nonlinear system does *not* have such a solution, we deduce that no such values $a, b, c, d$ and $e$ exist, hence no such differential equation is satisfied by $p_n(x)$, implying that the system is not a linear transformation of a classical orthogonal polynomial system.
Hence the whole question boils down to decide whether the given nonlinear system has non-trivial solutions, and to find these solutions in the affirmative case. As a matter of fact, with Gröbner bases methods, this question *can* be decided algorithmically [15]–[17]. Such an algorithm is implemented, e.g., in the computer algebra system REDUCE [16], and Maple's `solve` command can also solve such a system.
Note that the solution of the nonlinear system is not necessarily unique. For example, the Chebyshev polynomials of the first and second kind $T_n(x)$ and $U_n(x)$ satisfy the same recurrence equation, but a different differential equation. We will consider this example in more detail later.
If we apply this algorithm to the recurrence equation $p_{n+2}(x) - x\,p_{n+1}(x)$ of the power $p_n(x) = x^n$, it generates the complete solution set, given by Table 1:1 and 1:4a.
The following statement summarizes the above considerations.

**Algorithm 1** This algorithm decides whether a given holonomic three-term recurrence equation has classical orthogonal polynomial solutions, and returns their data if applicable. The algorithm is applicable to all entries of Table 1 independently of the orthogonality of the system under consideration.

1. **Input:** a holonomic three-term recurrence equation
$$q_n(x)\,p_{n+2}(x) + r_n(x)\,p_{n+1}(x) + s_n(x)\,p_n(x) = 0 \qquad (q_n(x), r_n(x), s_n(x) \in \mathbb{Q}[n,x]) \ .$$
2. **Shift:** Shift by $\max\{n \in \mathbb{N}_0 \,|\, n \text{ is zero of either } q_{n-1}(x) \text{ or } s_n(x)\} + 1$ if necessary.
3. **Rewriting:** Rewrite the recurrence equation in the form
$$p_{n+1}(x) = t_n(x)p_n(x) + u_n(x)\,p_{n-1}(x) \qquad (t_n(x), u_n(x) \in \mathbb{Q}(n,x)) \ .$$
If either $t_n(x)$ is not a polynomial of degree one in $x$ or $u_n(x)$ is not constant with respect to $x$, then return `"no classical orthogonal polynomial solution exists"`; exit.
4. **Standardization:** Given now $A_n, B_n$ and $C_n$ by
$$p_{n+1}(x) = (A_n\,x + B_n)\,p_n(x) - C_n\,p_{n-1}(x) \qquad (A_n, B_n, C_n \in \mathbb{Q}(n),\ A_n \neq 0) \ ,$$

define

$$\frac{k_{n+1}}{k_n} := A_n = \frac{v_n}{w_n} \qquad (v_n, w_n \in \mathbb{Q}[n])$$

according to (40).

5. **Make monic:** Set

$$\tilde{B}_n := \frac{B_n}{A_n} \in \mathbb{Q}(n) \qquad \text{and} \qquad \tilde{C}_n := \frac{C_n}{A_n\, A_{n-1}} \in \mathbb{Q}(n)$$

and bring them in lowest terms. If the degree of either the numerator or the denominator of $\tilde{B}_n$ is larger than 2, or if the degree of either the numerator or the denominator of $\tilde{C}_n$ is larger than 4, return `"no classical orthogonal polynomial solution exists"`; exit.

6. **Polynomial Identities:** Set

$$\tilde{B}_n = \frac{2bn\,(a(n-1)+d)+e(d-2a)}{(2a(n-1)+d)\,(2an+d)}$$

and

$$\tilde{C}_n = \frac{-n\,(a(n-2)+d)}{(a\,(2n-1)+d)(a\,(2n-3)+d)}\left(c+\frac{b(n-1)+e}{(2a\,(n-1)+d)^2}\Big((ae-bd)-ab\,(n-1)\Big)\right),$$

using the unknowns $a, b, c, d$ and $e$. Multiply these identities by their common denominators, and bring them therefore in polynomial form.

7. **Equating Coefficients:** Equate the coefficients of the powers of $n$ in the two resulting equations. This results in a nonlinear system in the unknowns $a, b, c, d$ and $e$. Solve this system by Gröbner bases methods. If the system has no solution, then return `"no classical orthogonal polynomial solution exists"`; exit.

8. **Output:** Return the classical orthogonal polynomial solutions of the differential equations (17) given by the solution vectors $(a, b, c, d, e)$ of the last step, according to the classification of Table 1, together with the information about the standardization given by (40). This information includes the density

$$\frac{\rho(x)}{C} = \frac{1}{\sigma(x)}\, e^{\int \frac{\tau(x)}{\sigma(x)}\, dx}$$

given by (11), and the interval by the zeros of $\sigma(x)$. □

We would like to give the following comments on the above algorithm:

1. The use of Gröbner bases is not always necessary. The following observation yields an ad hoc method to solve the nonlinear system. Observe that the coefficients of the powers of $n$ of the polynomial identity concerning $\tilde{B}_n$ of step 5 of the algorithm can be written using the variables

$$\{a^2, a\,d, a\,e, a\,b, d^2, d\,e, d\,b\} \tag{43}$$

Then all the derived equations are *linear* in the seven variables (43).

Furthermore the coefficients of the powers of $n$ of the polynomial identity concerning $\tilde{C}_n$ of step 5 of the algorithm essentially are products of *exactly* four terms out of $a, b, c, d, e$. Any of these can be written as a product of two of the variables

$$\{d_1 = a^2, d_2 = a\,d, d_3 = a\,e, d_4 = a\,b, d_5 = d^2, d_6 = d\,e, d_7 = d\,b, d_8 = a\,c, d_9 = d\,c\} \; . \tag{44}$$

This is the set of variables (43) plus the two variables $d_8, d_9$. All these equations are of *second order* in the nine variables (44).

Obviously the resulting system can be solved by first finding the solution space corresponding to the linear subsystem, which then can be substituted in the second order subsystem. The resulting second order system can be solved by ad hoc elimination (and possibly rational factorization).

If one has found the variables given by (44), then it is easy to calculate $a, b, c, d$ and $e$, or one realizes that no such solution exists.

2. Note moreover that, although Gröbner bases techniques apply only rational arithmetic, hence give rational solutions only, the technique described shows that solutions are also detected if they involve radicals.

3. If the input recurrence equation has further parameters, in step 6 of the algorithm one should solve for all variables including these additional ones, see Example 2.

**Example 1** As a first example, we consider the recurrence equation

$$(n+2)\,P_{n+2}(x) - x\,(n+1)\,P_{n+1}(x) + n\,P_n(x) = 0 \; .$$

Since $s_0(x) \equiv 0$, we see that the shift $p_n(x) := P_{n+1}(x)$ is necessary. For $p_n(x)$, we have the recurrence equation

$$(n+3)\,p_{n+2}(x) - x\,(n+2)\,p_{n+1}(x) + (n+1)\,p_n(x) = 0 \; . \tag{45}$$

In the first steps this recurrence equation is brought into the form

$$p_{n+1}(x) = \frac{n+1}{n+2}\,x\,p_n(x) - \frac{n}{n+2}\,p_{n-1}(x) \; ,$$

hence

$$A_n = \frac{k_{n+1}}{k_n} = \frac{n+1}{n+2} = \frac{v_n}{w_n} \; ,$$

and therefore

$$k_n = \frac{1}{n+1}\,k_0 \; .$$

Moreover, for monic $\tilde{p}_n(x) = p_n(x)/k_n$ we get

$$\tilde{p}_{n+1}(x) = x\,\tilde{p}_n(x) + \tilde{p}_{n-1}(x) \; ,$$

hence $\tilde{B}_n = 0$ and $\tilde{C}_n = 1$. In step 5 of the algorithm, the polynomial identity concerning $\tilde{B}_n$ then reads as

$$-2\,a\,b\,n^2 + (2\,a\,b - 2\,d\,b)\,n - d\,e + 2\,a\,e = 0 \; ,$$

leading to the linear system

$$a\,b=0,\quad 2\,a\,b-2\,d\,b=0,\quad d\,e-2\,a\,e=0$$

for the variables $a\,b, d\,b, d\,e, a\,e$, with the solution

$$\{d\,e=2\,a\,e, a\,b=0, d\,b=0\}\ , \tag{46}$$

$ae$ being arbitrary. After substituting the corresponding equations

$$\{d_6=2\,d_3, d_4=0, d_7=0\} \tag{47}$$

into the polynomial identity concerning $\tilde{C}_n$, we equate the coefficients and receive the second order equations

$$4\,d_1\ (4\,d_1+d_8)=0\ , \tag{48}$$

$$23\,d_1\,d_5-28\,d_1\,d_2+12\,{d_1}^2-8\,d_2\,d_5+{d_5}^2=0\ , \tag{49}$$

$$92\,{d_1}^2-96\,d_1\,d_2+24\,d_1\,d_5+5\,d_2\,d_9-20\,d_1\,d_9+20\,d_1\,d_8+{d_3}^2=0\ , \tag{50}$$

$$92\,d_1\,d_2-56\,{d_1}^2-48\,d_1\,d_5+8\,d_2\,d_5-6\,d_2\,d_9+d_5\,d_9+12\,d_1\,d_9-8\,d_1\,d_8=0\ , \tag{51}$$

$$-\,8\,d_1\ (8\,d_1-4\,d_2-d_9+2\,d_8)=0 \tag{52}$$

in terms of the variables (44). The first of these equations leads to two possibilities: either $d_1=0$ or $d_8=-4\,d_1$. One realizes quickly that the first of these possibilities implies $a=d=0$ which is not allowed. Hence we must have

$$d_8=-4\,d_1,\qquad \text{or}\qquad c=-4\,a\ , \tag{53}$$

and $d_1\neq 0$, i.e. $a\neq 0$. At this point we have already determined $\sigma(x)$ since by (46) one has $a\,b=0$, hence $b=0$ and therefore

$$\sigma(x)=ax^2+bx+c=a(x^2-4)\ .$$

Hence possible orthogonal polynomial solutions of (45) are defined in the interval $[-2,2]$.
We substitute now (53) in (49)–(52). Then the last equation reads as

$$-8\,d_1\ (-4\,d_2-d_9)=0\ .$$

Since $d_1\neq 0$, we conclude that

$$d_9=-4\,d_2\qquad \text{or}\qquad c=-4\,a\ . \tag{54}$$

In terms of $a,b,c,d$ and $e$ this yields nothing new, but it shows the compatibility of (52) with (48).
Substituting (54) in (48)–(52) gives two trivial identities, and three complicated ones. In these three equations, we finally resubstitute the original variables by (44), and after a rational factorization we get

$$(\,3\,a-d\,)\,(\,a-d\,)\,(\,2\,a-d\,)^2=0\ ,$$

$$-4\,a\,(\,3\,a-d\,)\,(\,a-d\,)\,(\,2\,a-d\,)=0\ ,$$

$$a^2\,\left(12\,a^2-16\,a\,d+4\,d^2+e^2\right)=0\ .$$

Hence either

$$d = a, \quad \text{or} \quad d = 2a, \quad \text{or} \quad d = 3a \; .$$

In the first of these cases, one gets $e = 0$ and the differential equation

$$(x^2 - 4)\, y''(x) + x\, y'(x) - n(n-2)\, y(x) = 0 \tag{55}$$

corresponding to the density

$$\rho(x) = -\frac{1}{\sigma(x)}\, e^{\int \frac{\tau(x)}{\sigma(x)}\, dx} = \frac{1}{\sqrt{4 - x^2}} \; .$$

The corresponding orthogonal polynomials are multiples of translated Chebyshev polynomials of the first kind

$$p_n(x) = k_n\, C_n(x) = \frac{p_0}{n+1}\, C_n(x) = \frac{2\, p_0}{n+1}\, T_n(x/2) \qquad (n \geq 0) \tag{56}$$

(see e.g. [1], Table 22.2, and (22.5.11); $C_n(x)$ are monic, but $C_0 = 2$, see also Table 22.7), hence finally

$$P_n(x) = p_{n-1}(x) = \frac{2\, P_1}{n}\, T_{n-1}(x/2) \qquad (n \geq 1) \; .$$

In the second of the above cases, i.e. for $d = 2a$, one gets the equation

$$a^2\, (e - 2a)\, (e + 2a) = 0$$

with two possible solutions $e = \pm 2a$ that give the differential equations

$$(x^2 - 4)\, y''(x) + 2\, (x+1)\, y'(x) - n(n-3)\, y(x) = 0 \; , \tag{57}$$

and

$$(x^2 - 4)\, y''(x) + 2\, (x-1)\, y'(x) - n(n-3)\, y(x) = 0 \; . \tag{58}$$

They correspond to the densities

$$\rho(x) = \sqrt{\frac{4+x}{4-x}} \qquad \text{and} \qquad \rho(x) = \sqrt{\frac{4-x}{4+x}} \; ,$$

respectively, hence the orthogonal polynomials are multiples of the Jacobi polynomials $P_n^{(1/2,-1/2)}(x/2)$ and $P_n^{(-1/2,1/2)}(x/2)$.
Finally, in the third of the above cases, i.e. for $d = 3a$, we get again $e = 0$ and

$$(x^2 - 4)\, y''(x) + 3\, x\, y'(x) - n(n-4)\, y(x) = 0 \tag{59}$$

corresponding to the density

$$\rho(x) = -\frac{1}{\sigma(x)}\, e^{\int \frac{\tau(x)}{\sigma(x)}\, dx} = \sqrt{4 - x^2} \; .$$

The corresponding orthogonal polynomials are multiples of translated Chebyshev polynomials of the second kind

$$p_n(x) = k_n\, S_n(x) = \frac{p_0}{n+1}\, S_n(x) = \frac{p_0}{n+1}\, U_n(x/2) \qquad (n \geq 0) \tag{60}$$

(see e.g. [1], Table 22.2, and (22.5.13); $S_n(x)$ are monic, see also Table 22.8), hence

$$P_n(x) = p_{n-1}(x) = \frac{P_1}{n} U_{n-1}(x/2) \qquad (n \geq 1) \; .$$

We see that the recurrence equation (45) has four different classical orthogonal polynomial solutions!

**Example 2** As a second example, we consider the recurrence equation

$$p_{n+2}(x) - (x - n - 1)\, p_{n+1}(x) + \alpha\, (n+1)^2\, p_n(x) = 0 \tag{61}$$

depending on the parameter $\alpha \in \mathbb{R}$. Here obviously the question arises whether or not there are any instances of this parameter for which there are classical orthogonal polynomial solutions. In step 6 of Algorithm 1 we therefore solve also for this unknown parameter. This gives a slightly more complicated nonlinear system, with the unique solution

$$\left\{ b = 2\,c, c = c, d = -4\,c, e = 0, a = 0, \alpha = \frac{1}{4} \right\} \; .$$

Hence the only possible value for $\alpha$ with classical orthogonal polynomial solutions is $\alpha = 1/4$, in which case one gets the differential equation

$$\left( x + \frac{1}{2} \right) p_n''(x) - 2\,x\,p_n'(x) - 2\,n\,p_n(x) = 0$$

with density

$$\rho(x) = 2\,e^{-2x}$$

in the interval $[-1/2, \infty]$, corresponding to shifted Laguerre polynomials.

# 4 Application: The Legendre Addition Theorem

As an application of Algorithm 1 in this section we show how the particular case

$$P_n(x^2 + (1-x^2)\cos\theta) = P_n(x)^2 + 2 \sum_{k=1}^{n} \frac{(n-k)!}{(n+k)!} P_n^k(x)^2 \, \cos k\theta \; . \tag{62}$$

of the Legendre addition theorem ([14], 5.4.4, p. 239, see also [22], [11])

$$P_n\left(xy + \sqrt{1-x^2}\sqrt{1-y^2}\,\cos\theta\right) = P_n(x)\,P_n(y) + 2 \sum_{k=1}^{n} \frac{(n-k)!}{(n+k)!} P_n^k(x)\,P_n^k(y)\,\cos k\theta$$

can be *deduced* by linear algebra techniques. Note that (62) played an essential role in Weinstein's proof of the Bieberbach conjecture [22]. Here $P_n(x) = P_n^{(0,0)}(x)$ are the *Legendre polynomials*, and $P_n^k(x)$ are called the *associated Legendre functions*. Our goal will be to identify these functions. In our deduction, we partially follow [6], see also the first author's review [9]. For the given purpose, we write

$$P_n(x^2 + (1-x^2)\cos\theta) = B_n^0(x) + 2 \sum_{k=1}^{n} \frac{(n-k)!}{(n+k)!} (1-x^2)^k\, B_n^k(x)\,\cos k\theta \tag{63}$$

with still unknown functions $B_n^k(x)$. Multiplying by $z^n$, and summing for $n = 0, 1, \ldots$ yields the generating function of the Legendre polynomials, hence ([1], (22.9.12))

$$\frac{1}{\sqrt{1 - z(2x^2 + (1 - x^2)(w+1/w)) + z^2}} = \sum_{n=0}^{\infty} \sum_{k=-n}^{n} \frac{(n-k)!}{(n+k)!} (1 - x^2)^k \, B_n^k(x) w^k z^n \tag{64}$$

where we put $w = e^{i\theta}$ and $B_n^{-k}(x) = B_n^k(x)$. In the sequel we consider this equation as a formal Laurent series expansion w.r.t. the variables $w$ and $z$. The functions $B_n^k(x)$ can be iteratively calculated by series approximations of the left hand function (e.g., using Maple), and it turns out that, for $0 \le k \le n \le 10$, for example, these form polynomials that are squares of another system of polynomials

$$D_n^k(x)^2 = B_n^k(x) \; . \tag{65}$$

We normalize $D_n^k(x)$ such that the highest coefficient has sign $(-1)^k$ (to be consistent with the definitions given in ([1], § 8)).
Now, we would like to find a three term recurrence equation w.r.t. $n$ valid for the polynomials $D_n^k(x)$.
For this purpose, we "guess" that

$$(ak + bn + c) \, D_{n+2}^k(x) + (dk + en + f) \, D_{n+1}^k(x) + (gk + hn + i) \, D_n^k(x) = 0$$

with unknowns $a, b, c, d, e, f, g, h, i$. Substituting the given values $D_n^k(x)$ $(0 \le k \le n \le 10)$ into this proposed recurrence equation yields a linear system which turns out to be consistent (although we have $\binom{9}{2} = 36$ equations, but only 9 unknowns), with the unique solution

$$(n - k + 2) \, D_{n+2}^k(x) - (2n + 3) \, x \, D_{n+1}^k(x) + (n + k + 1) \, D_n^k(x) = 0 \; . \tag{66}$$

Currently this recurrence equation is not yet proved, but this will be done soon. Assume for the moment that $E_n^k(x)$ are solutions of (66). Then by another application of linear algebra (see e.g. [20], [12]), this recurrence equation can be "squared", i.e. it is possible to calculate the recurrence equation of third order valid for $E_n^k(x)^2$. This step can be accomplished, e.g., by the procedure `rec*rec` of the `gfun` packacke ([20], see also [10]) with Maple, and results in the recurrence equation

$$\begin{aligned}
&(2n+5)(k+n+2)(k+n+1)^2 \, S_n^k(x) \\
&\quad - (2n+3)(k+n+2)(k^2 - n^2 + 4x^2n^2 - 4n + 16x^2n + 15x^2 - 4) \, S_{n+1}^k(x) \\
&\quad - (2n+5)(-2+k-n)(k^2 - n^2 + 4x^2n^2 - 4n + 16x^2n + 15x^2 - 4) \, S_{n+2}^k(x) \\
&\quad + (2n+3)(-2+k-n)(-3+k-n)^2 \, S_{n+3}^k(x) = 0
\end{aligned} \tag{67}$$

for the squares $S_n^k(x) = E_n^k(x)^2$.
If we are able to prove that this recurrence equation is valid for our unknown functions $B_n^k(x)$, then, by an a posteriori argument, we have deduced (65), since we have luckily found the "square root" recurrence equation (66) of (67).
Next we show how it can be discovered independently that $B_n^k(x)$ satisfy (67). We can rewrite (64) as

$$F_n^k(x) = \frac{(n-k)!}{(n+k)!} (1 - x^2)^k \, B_n^k(x) = \mathrm{CT}_{z,w} \, G_n^k(z, w)$$

with

$$G_n^k(z,w) := \frac{1}{\sqrt{1-z(2x^2+(1-x^2)(w+1/w))+z^2}}\,\frac{1}{z^n\,w^k}$$

where $\mathrm{CT}_{z,w}\,G_n^k(z,w)$ denotes the constant term of the double Laurent series $G_n^k(z,w)$. To obtain a recurrence equation for $F_n^k(x)$, we try to find polynomials $p_0, p_1, p_2, p_3$ in the variables $n, k$, and $x$, and polynomials $G_1$, and $G_2$, both of degree 2 in both $z$ and $w$, such that

$$p_0\,G_n^k + p_1\,G_{n+1}^k + p_2\,G_{n+2}^k + p_3\,G_{n+3}^k - z\,\frac{\partial}{\partial z}\left(\frac{G_1\,G_n^k}{z^3\,w}\right) - w\,\frac{\partial}{\partial w}\left(\frac{(1+z)\,G_2\,G_n^k}{z^3\,w}\right) = 0\;.$$

Substituting $G_1$ and $G_2$ generically, and dividing by $G_n^k$, a polynomial identity is derived, and by equating coefficients w.r.t. $z$ and $w$, we get a linear system, again. Solving this system results in the identity

$$\begin{aligned}
&-(\,n+3\,)\,(\,k+n+2\,)\,(\,2+n-k\,)\,G_n^k(z,w)\\
&\quad+(\,n+2\,)\,(\,4\,n^2\,x^2-n^2+22\,n\,x^2-6\,n-9+30\,x^2+k^2\,)\,G_{n+1}^k(z,w)\\
&\quad-(\,4\,n^2\,x^2-n^2+18\,n\,x^2-4\,n+20\,x^2+k^2-4\,)\,(\,n+3\,)\,G_{n+2}^k(z,w)\\
&\quad+(\,n+2\,)\,(\,n+3+k\,)\,(\,n+3-k\,)\,G_{n+3}^k(z,w) =\\
&\quad z\,\frac{\partial}{\partial z}\left(\frac{G_1(z,w)\,G_n^k(z,w)}{z^3\,w}\right)+w\,\frac{\partial}{\partial w}\left(\frac{(1+z)\,G_2(z,w)\,G_n^k(z,w)}{z^3\,w}\right)
\end{aligned}\tag{68}$$

for certain polynomials $G_1(z,w)$ and $G_2(z,w)$ which are reproduced in the appendix.
Since *any* formal Laurent series $f(z)$ satisfies

$$z\,\frac{\partial}{\partial z}f(z) = 0\;,$$

applying $\mathrm{CT}_{z,w}$ to the identity (68) yields the three term recurrence equation

$$\begin{aligned}
&-(\,n+3\,)\,(\,k+n+2\,)\,(\,2+n-k\,)\,F_n^k(x)\\
&\quad+(\,n+2\,)\,(\,4\,n^2\,x^2-n^2+22\,n\,x^2-6\,n-9+30\,x^2+k^2\,)\,F_{n+1}^k(x)\\
&\quad-(\,4\,n^2\,x^2-n^2+18\,n\,x^2-4\,n+20\,x^2+k^2-4\,)\,(\,n+3\,)\,F_{n+2}^k(x)\\
&\quad+(\,n+2\,)\,(\,n+3+k\,)\,(\,n+3-k\,)\,F_{n+3}^k(x) = 0
\end{aligned}$$

for $F_n^k(x)$. Using the recurrence equation

$$(1+k+n)\,a_n + (-1+k-n)\,a_{n+1} = 0$$

which is valid for the factor $(n+k)!/(n-k)!$, one can finally show (e.g., with `rec*rec` again) that $B_n^k(x)$ satisfy the recurrence equation (67), as announced.
Now we know that $B_n^k(x)$ are the squares of the polynomial system $D_n^k(x)$ defined by (66). But who are these polynomials? A direct application of Algorithm 1 is not possible since our polynomials $D_n^k(x)$ live for $n \geq k$ rather than for $n \geq 0$. Hence, we deal with

$$T_n^k(x) := D_{n+k}^k(x)$$

for which, by (66), we get the recurrence equation

$$(n+2)\,T_{n+2}^k(x) - (2n+2k+3)\,x\,T_{n+1}^k(x) + (n+2k+1)\,T_n^k(x) = 0\;.$$

Now an application of Algorithm 1 shows that the differential equation

$$(x^2-1)\,T_n^{k''}(x) + 2\,(1+k)\,x\,T_n^{k'}(x) - n\,(n-2k-3)\,T_n^k(x) = 0$$

is valid, and that $\frac{k_{n+1}}{k_n} = \frac{2n+2k+1}{n+1}$. Hence one has the density

$$\rho(x) = (1-x^2)^k \ ,$$

i.e., we have multiples of Jacobi-Gegenbauer polynomials

$$T_n^k(x) = \kappa(k)\,\frac{2^n\,(k+1/2)_n}{n!}\,Q_n^{(k,k)}(x)$$

with some function $\kappa(k)$ not depending on $n$, and

$$Q_n^{(\alpha,\beta)}(x) = \frac{1}{\binom{2n+\alpha+\beta}{n}} \sum_{k=0}^{n} \binom{n+\alpha}{k}\binom{n+\beta}{n-k}\,(x+1)^k\,(x-1)^{n-k}$$

denoting the monic Jacobi polynomials (see e.g. [1], (22.3.1)).
To determine $\kappa(k)$, we compare the coefficients of $\cos^n\theta$ in (63). Since ([1], (22.3))

$$P_n(x) = \frac{(2n)!}{2^n\,n!^2}\,x^n + O(x^{n-1}) \qquad \text{and} \qquad T_n(x) = 2^{n-1}\,x^n + O(x^{n-1}) \ ,$$

we get

$$\begin{aligned} P_n(x^2+(1-x^2)\cos\theta) &= \frac{(2n)!}{2^n\,n!^2}\Big(x^2+(1-x^2)\cos\theta\Big)^n + O(\cos^{n-1}(\theta)) \\ &= \frac{(2n)!}{2^n\,n!^2}\,(1-x^2)^n\,\cos^n\theta + O(\cos^{n-1}(\theta)) \\ &= 2\,\frac{1}{(2n)!}\,(1-x^2)^n\,D_n^n(x)^2\,\cos n\theta \\ &= \frac{2}{(2n)!}\,(1-x^2)^n\,D_n^n(x)^2\,2^{n-1}\cos^n\theta + O(\cos^{n-1}(\theta)) \ , \end{aligned}$$

because $\cos n\theta = T_n(\cos\theta)$. Hence, we arrive at

$$D_n^n(x)^2 = \frac{(2n)!^2}{2^{2n}\,n!^2}\,(1-x^2)^n \ ,$$

and, setting $n=k$, we see that we must choose $\kappa(k) = (-1)^k\,\frac{(2k)!}{2^k\,k!}$, according to the normalization with highest coefficient sign $(-1)^k$ of $D_n^k(x)$. We have finally determined the representation

$$D_n^k(x) = T_{n-k}^k(x) = (-1)^k\,\frac{(2k)!}{2^k\,k!}\,\frac{2^{n-k}\,(k+1/2)_{n-k}}{(n-k)!}\,Q_{n-k}^{(k,k)}(x) \ ,$$

and therefore

$$P_n^k(x) = (-1)^k\,(1-x^2)^{k/2}\,\frac{2^{n-2k}\,(2k)!}{k!}\,\frac{(k+1/2)_{n-k}}{(n-k)!}\,Q_{n-k}^{(k,k)}(x) \ .$$

This finishes the proposed identification of $P_n^k(x)$.

# 5 Classical Orthogonal Polynomials of a Discrete Variable

In this section, we give similar results for classical orthogonal polynomials of a discrete variable, see Chapter 2 of [19]. These are given by a difference equation

$$\sigma(x)\,\Delta\nabla y(x) + \tau(x)\,\Delta y(x) + \lambda_n\, y(x) = 0\ , \tag{69}$$

where

$$\Delta y(x) = y(x+1) - y(x) \qquad \text{and} \qquad \nabla y(x) = y(x) - y(x-1)$$

denote the forward and backward difference operators, respectively. The classical discrete orthogonal polynomials are the polynomial solutions of the difference equation (69).
These polynomials can be classified similarly as in the continuous case according to the functions $\sigma(x)$ and $\tau(x)$; up to linear transformations the classical discrete orthogonal polynomials are classified acccording to Table 2 (see [19], Chapter 2, and § 6). We added the other possible solutions to those given in ([19], Chapter 2). In particular, case (2a) corresponds to the non-orthogonal solution $x^n$ in Table 1. Similarly as for the powers

$$\frac{d}{dx}x^n = n\,x^{n-1}\ ,$$

the *falling factorials* $x^{\underline{n}} := x(x-1)\cdots(x-n+1)$ satisfy

$$\Delta x^{\underline{n}} = n\,x^{\underline{n-1}}\ .$$

It turns out that they are connected with the Charlier polynomials by the limiting process

$$\lim_{\mu\to 0}(-1)^n\,\mu^n\,c_n^{(\mu)}(x) = \lim_{\mu\to 0}(x-n+1)_n\;{}_1F_1\left(\left.\begin{array}{c}-n\\ x-n+1\end{array}\right|\mu\right) = (-1)^n\,(x-n+1)_n = x^{\underline{n}}$$

where we used the hypergeometric representation given in ([19], (2.7.9)).
Note, however, that other than in the differential equation case the above type of difference equation is *not* invariant under general linear transformations, but only under integer shifts. We will have to take this under consideration.

The classical discrete orthogonal polynomial systems correspond to a discrete weight function $\rho(x)$, with ([19], Equation (2.4.1))

$$\frac{\rho(x+1)}{\rho(x)} = \frac{\sigma(x)+\tau(x)}{\sigma(x+1)}\ , \tag{70}$$

given by the Pearson type difference equation ([19], Equation (2.1.17))

$$\Delta\Big(\sigma(x)\,\rho(x)\Big) = \tau(x)\,\rho(x)\ . \tag{71}$$

As in the continuous case, by Lemma 1 (the measure $\mu$ is discrete now) they satisfy a recurrence equation (4)

$$p_{n+1}(x) = (A_n\,x + B_n)\,p_n(x) - C_n\,p_{n-1}(x)\ .$$

Multiplying (69) by $\rho(x)$, by (71) the difference equation takes the selfadjoint form ([19], Equation (2.1.18))

$$\Delta\Big(\sigma(x)\,\rho(x)\,\nabla y(x)\Big) + \lambda_n\,y(x) = 0\ ,$$

1. $\sigma(x) = 1,\ \tau(x) = \alpha\, x + \beta \quad \Longrightarrow \quad p_n(x) = K_n^{(\alpha,\beta)}(x)$, see § 6,

2a. $\sigma(x) = x,\ \sigma(x) + \tau(x) = 0 \quad \Longrightarrow \quad p_n(x) = x^{\underline{n}} := x(x-1)\cdots(x-n+1)$,

2b. $\sigma(x) = x,\ \sigma(x) + \tau(x) = \mu\ (\mu \neq 0) \quad \Longrightarrow \quad p_n(x) = c_n^{(\mu)}(x)$, the *Charlier polynomials*,

3 $\sigma(x) = x,\ \sigma(x)+\tau(x) = \mu\,(\gamma+x) \quad \Longrightarrow \quad p_n(x) = m_n^{(\gamma,\mu)}(x)$, the *Meixner polynomials*,

4 $\sigma(x) = x,\ \sigma(x)+\tau(x) = \frac{p}{1-p}\,(N-x) \ \Rightarrow\ p_n(x) = k_n^{(p)}(x,N)$, the *Krawchouk polynomials*,

5 $\sigma(x) = x(N+\alpha-x),\ \ \sigma(x)+\tau(x) = (x+\beta+1)(N-1-x) \ \Rightarrow\ p_n(x) = h_n^{(\alpha,\beta)}(x,N)$, the *Hahn polynomials (of first kind)*,

6 $\sigma(x) = x(x+\mu),\ \ \sigma(x)+\tau(x) = (\nu+N-1-x)(N-1-x) \ \Rightarrow\ p_n(x) = \tilde{h}_n^{(\mu,\nu)}(x,N)$, the *Hahn polynomials (of second kind)*.

Table 2: Normal Forms of Discrete Polynomials

from which one can deduce (as in Lemma 2) that $p_n(x)$ satisfy a difference rule (see also [19], (2.2.10))

$$\sigma(x)\,\nabla p_n(x) = \alpha_n\, p_{n+1}(x) + \beta_n\, p_n(x) + \gamma_n\, p_{n-1}(x)\ . \tag{72}$$

Analogous to the continuous case, one has a *Rodrigues formula* ([19], (2.2.8))

$$p_n(x) = \frac{E_n}{\rho(x)}\,\Delta^n\left(\rho(x)\prod_{k=0}^{n-1}\sigma(x-k)\right) \tag{73}$$

for some $E_n$ independent of $x$, and the *orthogonality relations* read as ([19], (2.3.1))

$$\int\limits_{\mathbb{R}} p_n(x)\,p_m(x)\,d\mu(x) = \sum_{x\in X} p_n(x)\,p_m(x)\,\rho(x) = \begin{cases} 0 & \text{if } m \neq n \\ h_n \neq 0 & \text{if } m = n \end{cases}$$

for some discrete set $X \subset \mathbb{R}$.
To deduce formulas similar to those in Theorem 1, we substitute (25)

$$p_n(x) = k_n\, x^n + k_n'\, x^{n-1} + k_n''\, x^{n-2} + \ldots$$

in the difference equations for $p_n(x)$ and for $p_{n+1}(x)$, in the recurrence equation and in the difference rule. Equating the three highest coefficients in any of these four equations yields twelve nonlinear equations in the twelve unknowns

$$A_n, B_n, C_n, \alpha_n, \beta_n, \gamma_n, \lambda_n, \lambda_{n+1}, k_n', k_{n+1}', k_n'', k_{n+1}''\ .$$

In particular the highest coefficient of the difference equation yields

$$\lambda_n = -(an(n-1) + dn)$$

again. We can assume that $\lambda_n \neq 0$ for $n \in \mathbb{N}$, hence $a(n-1) + d \neq 0$ for $n \in \mathbb{N}$ since otherwise no orthogonal polynomial solutions can exist; in particular, $d \neq 0$.
By a Gröbner basis calculation (invoked by the `solve` command of Maple or REDUCE) it turns out that there is a unique solution of the above nonlinear system, given by **Theorem 2** Let $p_n(x) = k_n\, x^n + \ldots$ $(n \in \mathbb{N}_0)$ be a family of polynomial solutions of the system of difference equations (69) that are orthogonal with respect to a discrete weight function $\rho(x)$. Then the difference rule (72)

$$\sigma(x)\, \nabla p_n(x) = \alpha_n\, p_{n+1}(x) + \beta_n\, p_n(x) + \gamma_n\, p_{n-1}(x)$$

is valid with

$$\alpha_n = a\, n\, \frac{k_n}{k_{n+1}}$$

$$\beta_n = -\frac{n\,(a(n-1)+d)\,(2(n-1)\,a\,(a(n-1)+a+d)+ad+2ae-bd)}{(2\,an+d)\,(2\,a(n-1)+d)} ,$$

$$\gamma_n = \frac{k_n}{k_{n-1}} \cdot \frac{n\,\,(a(n-1)+d)\,(a(n-2)+d)}{(a(2n-3)+d)\,(a(2n-1)+d)\,(2\,a(n-1)+d)^2}$$
$$\cdot \Big((n-1)(a(n-1)+d)(a^2(n-1)^2+a(n-1)d+4ac+2ae-bd-b^2)+ae^2-bde+cd^2\Big) ,$$

and the recurrence equation (4)

$$p_{n+1}(x) = (A_n\, x + B_n)\, p_n(x) - C_n\, p_{n-1}(x)$$

is valid with

$$A_n = \frac{k_{n+1}}{k_n} , \tag{74}$$

$$B_n = \frac{k_{n+1}}{k_n} \cdot \frac{(n-1)(d+2b)(a(n-1)+a+d)-2ae+d^2+de+2bd}{(2\,an+d)\,(2\,a(n-1)+d)} \tag{75}$$

and

$$C_n = -\frac{k_{n+1}}{k_n}\, \frac{\gamma_n}{a(n-1)+d} . \qquad \Box \tag{76}$$

Note that with an immense effort one can deduce Theorem 2 also by hand calculations, using a similar technique as in our proof of Theorem 1.
Theorem 2 is also valid for case (1) of Table 2, see § 6; the recurrence equation part is also valid in case (2a) of Table 2, and generates the recurrence equation $p_{n+2}(x) - (x-n-1)\, p_{n+1}(x) = 0$ for the falling factorial $p_n(x) = x^{\underline{n}}$.

As a byproduct of Theorem 2, we get for the ratio $k'_{n+1}/k'_n$

**Corollary 6** Let $p_n(x) = k_n\, x^n + k'_n\, x^{n-1} + \ldots$ $(n \in \mathbb{N}_0)$ be a family of orthogonal polynomial solutions of the system of difference equations (69). Then the relation

$$\frac{k'_{n+1}}{k'_n} = \frac{n+1}{n}\, \frac{(\,d\,n + 2\,b\,n + 2\,e\,)\,(\,2\,a\,(n-1) + d\,)}{(\,d + 2\,a\,n\,)\,(\,2\,b\,(n-1) + d\,n + 2\,e - d\,)}\, \frac{k_{n+1}}{k_n}$$

is valid. $\Box$

As in the continuous case, the Rodrigues constant $E_n$, given by (73), can be determined.

**Corollary 7** Let $p_n(x) = k_n\,x^n + \ldots$ $(n \in \mathbb{N}_0)$ be a family of orthogonal polynomial solutions of the system of difference equations (69) corresponding to a discrete weight function $\rho(x)$. Then $p_n(x)$ have a Rodrigues representation (73), $E_n$ being a hypergeometric term, satifying the relation

$$\frac{E_{n+1}}{E_n} = \frac{(a(n-1)+d)}{(a(2n-1)+d)\,d} \cdot \frac{k_{n+1}}{k_n}\;. \tag{77}$$

*Proof:* In ([19], (2.2.10)) it was shown that

$$\sigma(x)\,\nabla p_n(x) = \frac{\lambda_n}{n\,\tau_n'(0)}\left(\tau_n(x)\,p_n(x) - \frac{E_n}{E_{n+1}}\,p_{n+1}\right) \tag{78}$$

where $\tau_n(x) = \tau(x+n) + \sigma(x+n) - \sigma(x)$ ([19], (2.1.15)). An application of a method similar to that used in Theorem 2 yields the coefficients for this type of difference rule, and hence determines $E_{n+1}/E_n$ according to (77). $\square$

Also analogously we get

**Corollary 8** Let $p_n(x) = k_n\,x^n + \ldots$ $(n \in \mathbb{N}_0)$ be a family of orthogonal polynomial solutions of the system of difference equations (69) corresponding to a discrete weight function $\rho(x)$. Then the relation

$$\begin{aligned}\frac{h_{n+1}}{h_n} &= (\,d + a\,n - a\,)(-a^3\,n^4 - 4\,d\,c\,n\,a + d\,b\,e + d\,b^2\,n - d^2\,c - a\,e^2 - 4\,c\,n^2\,a^2 \\ &\quad - 2\,n\,a\,d\,e - 2\,n^3\,a^2\,d + n^2\,a\,d\,b + b\,n\,d^2 - n^2\,a\,d^2 + n^2\,b^2\,a - 2\,a^2\,e\,n^2)(\,n+1\,) \\ &\quad \Big/ \Big((\,d + 2\,a\,n + a\,)\,(\,d + 2\,a\,n - a\,)\,(\,d + 2\,a\,n\,)^2\Big) \cdot \left(\frac{k_{n+1}}{k_n}\right)^2\end{aligned}$$

is valid.

*Proof:* In ([19], (2.5.6)) it was proved that

$$\frac{C_n}{A_n} = \frac{k_{n-1}}{k_n} \cdot \frac{h_n}{h_{n-1}}\;,$$

hence

$$\frac{h_{n+1}}{h_n} = \frac{C_{n+1}}{A_{n+1}} \cdot \frac{k_{n+1}}{k_n}\;.$$

An application of Theorem 2 yields the result. $\square$

Similarly as in the continuous case, Theorem 2 can be used to generate an algorithm to test whether or not a given holonomic recurrence equation has classical discrete orthogonal polynomial solutions.

**Algorithm 2** This algorithm decides whether a given holonomic three-term recurrence equation has classical discrete orthogonal polynomial solutions, and returns their data if applicable.

1. **Input:** a holonomic three-term recurrence equation

   $$q_n(x)\,p_{n+2}(x) + r_n(x)\,p_{n+1}(x) + s_n(x)\,p_n(x) = 0 \qquad (q_n(x), r_n(x), s_n(x) \in \mathbb{Q}[n,x])\;.$$

2. **Shift:** Shift by $\max\{n \in \mathbb{N}_0\,|\,n$ is zero of either $q_{n-1}(x)$ or $s_n(x)\} + 1$ if necessary.

3. **Rewriting:** Rewrite the recurrence equation in the form

$$p_{n+1}(x) = t_n(x) p_n(x) + u_n(x)\, p_{n-1}(x) \qquad (t_n(x), u_n(x) \in \mathbb{Q}(n,x)) \ .$$

If either $t_n(x)$ is not a polynomial of degree one in $x$ or $u_n(x)$ is not constant with respect to $x$, return `"no classical discrete orthogonal polynomial solution exists"`; exit.

4. **Linear Transformation:** Rewrite the recurrence equation by the linear transformation $x \mapsto \frac{x-g}{f}$ with (yet) unknowns $f$ and $g$.

5. **Standardization:** Given now $A_n, B_n$ and $C_n$ by

$$p_{n+1}(x) = (A_n\, x + B_n)\, p_n(x) - C_n\, p_{n-1}(x) \qquad (A_n, B_n, C_n \in \mathbb{Q}(n),\ A_n \neq 0) \ ,$$

define

$$\frac{k_{n+1}}{k_n} := A_n = \frac{v_n}{w_n} \qquad (v_n, w_n \in \mathbb{Q}[n])$$

according to (74).

6. **Make monic:** Set

$$\tilde{B}_n := \frac{B_n}{A_n} \in \mathbb{Q}(n) \qquad \text{and} \qquad \tilde{C}_n := \frac{C_n}{A_n\, A_{n-1}} \in \mathbb{Q}(n)$$

and bring them in lowest terms. If the degree of either the numerator or the denominator of $\tilde{B}_n$ is larger than 2, if the degree of the numerator of $\tilde{C}_n$ is larger than 6, or if the degree of the denominator of $\tilde{C}_n$ is larger than 4, then return `"no classical discrete orthogonal polynomial solution exists"`; exit.

7. **Polynomial Identities:** Set

$$\tilde{B}_n = \frac{k_n}{k_{n+1}}\, B_n$$

according to (75), and

$$\tilde{C}_n = \frac{k_{n-1}}{k_{n+1}}\, C_n$$

according to (76), in terms of the unknowns $a, b, c, d, e, f$ and $g$. Multiply these identities by their common denominators, and bring them therefore in polynomial form.

8. **Equating Coefficients:** Equate the coefficients of the powers of $n$ in the two resulting equations. This results in a nonlinear system in the unknowns $a, b, c, d, e, f$ and $g$. Solve this system by Gröbner bases methods. If the system has no solution, then return `"no classical discrete orthogonal polynomial solution exists"`; exit.

9. **Output:** Return the classical orthogonal polynomial solutions of the difference equations (69) given by the solution vectors $(a, b, c, d, e, f, g)$ of the last step, according to the classification given in Table 2, together with the information about the standardization given by (74). This information includes the necessary linear transformation $fx + g$, as well as the discrete weight function $\rho(x)$ given by (70)

$$\frac{\rho(x+1)}{\rho(x)} = \frac{\sigma(x) + \tau(x)}{\sigma(x+1)} \ .$$

*Proof:* The proof is an obvious modification of Algorithm 1. The only difference is that we have to take a possible linear transformation $fx + g$ into consideration since the difference equation (69) is not invariant under those transformations. This leads to step 3 of the algorithm. $\square$

Note that an application of Algorithm 2 to the recurrence equation $p_{n+2}(x) - (x - n - 1)$ $p_{n+1}(x) = 0$ which is valid for the falling factorial $p_n(x) = x^{\underline{n}}$, generates the difference equation $x\Delta\nabla p_n(x) - x\Delta p_n(x) + np_n(x) = 0$ of Table 2:2a.

**Example 3** We consider again the recurrence equation (61)

$$p_{n+2}(x) - (x - n - 1)\, p_{n+1}(x) + \alpha\,(n + 1)^2\, p_n(x) = 0$$

depending on the parameter $\alpha \in \mathbb{R}$. This time, we are interested in classical discrete orthogonal polynomial solutions.
According to step 3 of Algorithm 2, we rewrite (61) using the linear transformation $x \mapsto \frac{x-g}{f}$ with yet unknowns $f$ and $g$. Step 4 yields the standardization

$$\frac{k_{n+1}}{k_n} = -\frac{1}{f}\;.$$

In step 7, we solve the resulting nonlinear system for the variables $\{a, b, c, d, e, f, g, \alpha\}$, resulting in

$$\left\{a = 0, b = b, c = -\frac{b\,(-e+d+b\,)}{d}, d = d, e = e, f = -\frac{d+2\,b}{d}, g = -\frac{e}{d}, \alpha = \frac{b\,(\,d+b\,)}{(\,d+2\,b\,)^2}\right\}. \quad (79)$$

This is a rational representation of the solution. Since we assume $\alpha$ to be arbitrary, we solve the last equation for $b$. This yields

$$b = -\frac{d}{2}\left(1 \pm \frac{1}{\sqrt{1-4\alpha}}\right)\;,$$

which cannot be represented without radicals. Substituting this into (79) yields the solution

$$\left\{a = 0, b = -\frac{d}{2}\left(1 \pm \frac{1}{\sqrt{1-4\alpha}}\right), c = \frac{4\,\alpha\,e-e-2\,\alpha\,d}{2(1-4\,\alpha)} \pm \frac{e}{2}\,\frac{1}{\sqrt{1-4\alpha}}, f = \mp\frac{1}{\sqrt{1-4\alpha}}, g = -\frac{e}{d}\right\},$$

$d$ and $e$ being arbitrary. It turns out that for $1 - 4\,\alpha > 0$ this corresponds to Meixner or Krawchouk polynomials.

**Example 4** Here we want to discuss the possibility that a given recurrence equation might have several classical discrete orthogonal solutions. Whereas the recurrence equation of the Hahn polynomials $h_n^{(\alpha,\beta)}(x, N)$ has (besides several linear transformations) only this single classical discrete orthogonal solution, the case $\beta = -\alpha$ results in two essentially different solutions.
Here one has the recurrence equation

$$\begin{aligned}&(n + 2 + \alpha)\,(2 + n)\,(2\,n + 2)\,(n - N + 1)\,p_{n+2}(x) + (3 + 2\,n)\\&\quad(-6\,n\,\alpha - 2\,n^2\,\alpha - 4\,n^2\,x - 12\,n\,x + 2\,n^2\,N + 6\,n\,N + 4\,N - 4\,\alpha - 8\,x)\,p_{n+1}(x)\\&\quad - (1 + n)\,(n + 1 - \alpha)\,(2\,n + 4)\,(n + N + 2)\,p_n(x) = 0\;.\end{aligned}$$

An application of Algorithm 2 shows that this recurrence equation corresponds to the two different difference equations

$$x\,(-x+1-\alpha+N)\,\Delta\,\nabla p_n(x)+(-2\,x+N+\alpha\,N)\,\Delta p_n(x)\,)+n\,(n-3)\,p_n(x)=0$$

and

$$(x+\alpha)\,(-x+1+N)\,\Delta\nabla p_n(x)-(\,2\,x-N+2\,\alpha+\alpha\,N)\,\Delta p_n(x)+n\,(n-3)\,p_n(x)=0\;.$$

## 6 A New Polynomial System

In this section, we would like to present a new polynomial system satisfying the difference equation (69), but not coming from a discrete weight function, see Table 2:1. We set $\sigma(x)=1$, and $\tau(x)=\alpha x+\beta$ ($\alpha\neq 0$). This yields the polynomials

$$K_n^{(\alpha,\beta)}(x)=\left(x+\frac{1+\beta}{\alpha}\right)_n\cdot\alpha^n\cdot{}_1F_1\left(\left.\begin{array}{c}-n\\ 1-x-n-\frac{1+\beta}{\alpha}\end{array}\right|-\frac{1}{\alpha}\right)=(-1)^n\cdot{}_2F_0\left(\left.\begin{array}{c}-n,x+\frac{1+\beta}{\alpha}\\ -\end{array}\right|\alpha\right)$$

where we have used $k_n=\alpha^n$. Note that the two different hypergeometric representations come from ([19], (2.7.7)–(2.7.8)); they convert into one another by changing the direction of summation. The functions $K_n^{(\alpha,\beta)}(x)$ are not covered in ([19], Chapter 2) since

$$\frac{\rho(x+1)}{\rho(x)}=1+\alpha\,x+\beta\;,\qquad\text{i.e.}\qquad\rho(x)=\alpha^x\,\left(\frac{1+\beta}{\alpha}\right)_x$$

does not correspond to a valid weight function $\rho(x)$ over a suitable discrete real set.
The deductions in [19] show that (77) is valid, and we find that the Rodrigues constant is $E_n=1$, hence (73)

$$K_n^{(\alpha,\beta)}(x)=\frac{1}{\rho(x)}\,\Delta^n\rho(x)=\frac{\alpha^{-x}}{((1+\beta)/\alpha)_x}\,\Delta^n\left(\alpha^x\,\left(\frac{1+\beta}{\alpha}\right)_x\right)\;.$$

Applying Zeilberger's algorithm ([24], see also [7]) to one of the hypergeometric representations yields the recurrence equation

$$K_{n+2}^{(\alpha,\beta)}(x)-(n\,\alpha+x\,\alpha+\alpha+\beta)\,K_{n+1}^{(\alpha,\beta)}(x)-(1+n)\,\alpha\,K_n^{(\alpha,\beta)}(x)=0\qquad(80)$$

with the initial functions

$$K_{-1}^{(\alpha,\beta)}(x)\equiv 0\;,\qquad K_0^{(\alpha,\beta)}(x)=1\;.$$

The next three polynomials are

$$K_1^{(\alpha,\beta)}(x)=\alpha\,x+\beta\;,\qquad K_2^{(\alpha,\beta)}(x)=x^2\alpha^2+\alpha\;(2\,\beta+\alpha)\,x+\beta^2+\alpha+\alpha\,\beta\;,$$

$$K_3^{(\alpha,\beta)}(x)=\alpha^3x^3+3\,\alpha^2(\alpha+\beta)x^2+\alpha\Big(6\,\alpha\,\beta+2\,\alpha^2+3\,\beta^2+3\,\alpha\Big)x+\beta^3+2\,\alpha^2\beta+3\,\alpha\,\beta^2+3\,\alpha\,\beta+2\,\alpha^2.$$

According to (78), the difference rules

$$\nabla K_n^{(\alpha,\beta)}(x)=-(\alpha\,x+\alpha\,n+\beta)\,K_n^{(\alpha,\beta)}(x)+K_{n+1}^{(\alpha,\beta)}(x)=\alpha\,n\,K_{n-1}^{(\alpha,\beta)}(x)\qquad(81)$$

are valid which are consistent with Theorem 2. Another application of Zeilberger's algorithm yields furthermore the recurrence equation

$$(1+\beta+x\,\alpha+\alpha)\,K_n^{(\alpha,\beta)}(x+2)-(\beta+x\,\alpha+n\,\alpha+2+\alpha)K_n^{(\alpha,\beta)}(x+1)+K_n^{(\alpha,\beta)}(x)=0$$

with respect to $x$ which also can be obtained combining (80) and (81).
Note that the given polynomial system completes the classification of the polynomial solutions of (69) given in ([19], Chapter 2), resulting in Table 2. Since the given system satisfies a difference rule which is compatible with Theorem 2, Algorithm 2 can also recognize these solutions.

# 7 Coefficient Recurrence Equation

By equating the coefficients of the powers of $x$ in their differential equation (17) one gets for the continuous classical orthogonal polynomials

$$p_n(x)=\sum_{k=0}^{n} a_k^{(n)}\,x^k \tag{82}$$

the three term recurrence equation ($k=0,\ldots,n-1$)

$$(n-k)\,(a\,(n+k-1)+d)\,a_k^{(n)}=(k+1)\,\Big((k\,b+e)\,a_{k+1}^{(n)}+(k+2)\,c\,a_{k+2}^{(n)}\Big)\,, \qquad (a_{n+1}^{(n)}\equiv 0)$$

which was given in ([13], Eq. (3)).
In this section, we discuss the same question in the discrete case. Since the operators $\Delta$ and $\nabla$ behave nicely with falling factorials only, the solution is a little more difficult in this setting. Throughout this section we assume an expansion (82), i.e., we use the renamings $k_n=a_n^{(n)}$, $k_n'=a_{n-1}^{(n)}$, etc.
Let now $p_n(x)$ be a classical discrete system satisfying the difference equation (69)

$$\sigma(x)\,\Delta\nabla p_n(x)+\tau(x)\,\Delta p_n(x)+\lambda_n\,p_n(x)=0$$

with $\sigma(x)=ax^2+bx+c$, $\tau(x)=dx+e$ and $\lambda_n=-(an(n-1)+dn)$. Note that the difference equation is equivalent to

$$\Big(\sigma(x)+\tau(x)\Big)\,\Delta p_n(x)=\sigma(x)\,\nabla p_n(x)-\lambda_n\,p_n(x)\;, \tag{83}$$

hence the use of the double difference operator $\Delta\nabla$ can be omitted.
Taylor's theorem yields the expansions

$$\Delta p_n(x)=\sum_{k=1}^{n}\frac{p_n^{(k)}(x)}{k!} \qquad \text{and} \qquad \nabla p_n(x)=\sum_{k=1}^{n}(-1)^{k+1}\frac{p_n^{(k)}(x)}{k!}\;,$$

and by Leibniz's formula we receive the representations

$$\Delta p_n(x)=\sum_{k=0}^{n-1}\sum_{j=1}^{n-k}\binom{k+j}{j}\,a_{k+j}^{(n)}\,x^k \qquad \text{and} \qquad \nabla p_n(x)=\sum_{k=0}^{n-1}\sum_{j=1}^{n-k}\binom{k+j}{j}\,(-1)^{j+1}\,a_{k+j}^{(n)}\,x^k$$

for $\Delta p_n(x)$ and $\nabla p_n(x)$.

Substituting these representations in (83) and equating the coefficients of the powers of $x$ results finally in the recurrence equations

$$(2a(n-1)+d)\,a_{n-1}^{(n)} = n\left(e+\left(b+\frac{d}{2}\right)(n-1)\right)\,a_n^{(n)}\ , \tag{84}$$

and for $0 \le k \le n-2$

$$(n-k)\,(a\,(n+k-1)+d)\,a_k^{(n)} = \frac{k+1}{2}\Big((2b+d)\,k+2e\Big)\,a_{k+1}^{(n)} + \tag{85}$$

$$\sum_{j=2}^{n-k}\left(a\,(1+(-1)^j)\binom{k+j}{j+2}+\Big(b\,(1-(-1)^j)+d\Big)\binom{k+j}{j+1}+\Big(c\,(1+(-1)^j)+e\Big)\binom{k+j}{j}\right)a_{k+j}^{(n)}\ .$$

Note that from (84) and (85) one can deduce the discrete equivalents of (29) and (33).
As an example, we apply these formulas to the falling factorials $p_n(x) = x^{\underline{n}}$ which are solutions of the difference equation (Table 2:2a)

$$x\Delta\nabla p_n(x) - x\Delta p_n(x) + n p_n(x) = 0\ .$$

Here we have $a = c = e = 0$, $b = -d = 1$. Hence for the coefficients $s(n,k)$ of

$$x^{\underline{n}} = \sum_{k=0}^{n} s(n,k)\,x^k$$

that are called the *Stirling numbers of the first kind*, we get the relation

$$(n-k)\,s(n,k) = \sum_{j=1}^{n-k}(-1)^j\binom{k+j}{j+1}s(n,k+j)$$

which is due to Lagrange (see e.g. [5], p. 215, Theorem C).

# Appendix

Using higher coefficients $k_n'', k_n''', \ldots$ of the polynomials

$$y(x) = p_n(x) = k_n\,x^n + k_n'\,x^{n-1} + k_n''\,x^{n-2} + k_n'''\,x^{n-3} + k_n''''\,x^{n-4} + \ldots \qquad (n \in \mathbb{N}_0, k_n \neq 0)$$

that are solutions of (1)

$$\sigma(x)\,y''(x) + \tau(x)\,y'(x) + \lambda_n\,y(x) = 0\ ,$$

or (69)

$$\sigma(x)\,\Delta\nabla y(x) + \tau(x)\,\Delta y(x) + \lambda_n\,y(x) = 0\ ,$$

respectively, it is possible to use these as auxiliary variables, hence determining formulas for them. Note that to obtain the results of this article, it was crucial to make use of $k_n''$ which was never done before.
In this appendix, we collect results for the coefficients $k_n'', k_n'''$ and $k_n''''$ in both the continuous and discrete cases.

In the continuous case, we have the following properties:

$$\frac{k''_{n+1}}{k''_n} = \frac{(n^2 b^2 + 2 b e n + 2 c n a - b^2 n - b e + c d + e^2)}{(n^2 b^2 - 3 b^2 n + 2 b e n + 2 c n a + c d + 2 b^2 - 2 c a - 3 b e + e^2)} \cdot \frac{(n+1)(2 a n - 2 a + d)(d - 3 a + 2 a n)}{(d - a + 2 a n)(d + 2 a n)(n-1)} \frac{k_{n+1}}{k_n},$$

$$\frac{k'''_{n+1}}{k'''_n} = (n+1)(n^3 b^3 - 3 n^2 b^3 + 6 n^2 b c a + 3 n^2 b^2 e - 6 n b c a + 3 n b c d + 6 n e c a$$
$$- 6 n b^2 e + 3 n b e^2 + 2 n b^3 + 3 e c d - 2 b c d - 3 b e^2 - 2 c e a + 2 b^2 e + e^3)$$
$$(d - 3 a + 2 a n)(2 a n - 4 a + d) \Big/ \Big( (d + 2 a n)(d - a + 2 a n)(n-2)(n^3 b^3$$
$$+ 6 n^2 b c a + 3 n^2 b^2 e - 6 n^2 b^3 + 6 n e c a + 11 n b^3 - 18 n b c a - 12 n b^2 e$$
$$+ 3 n b e^2 + 3 n b c d - 6 b^3 - 5 b c d - 6 b e^2 + 11 b^2 e + 3 e c d - 8 c e a$$
$$+ 12 b c a + e^3) \Big) \frac{k_{n+1}}{k_n}$$

and

$$\frac{k''''_{n+1}}{k''''_n} = (n+1)(6 n^2 b^2 c d + 11 n^2 b^4 + 12 c^2 n d a + n^4 b^4 - 6 n^3 b^4 - 36 n^2 b^2 c a$$
$$+ 12 n c e^2 a + 22 n b^3 e - 18 n b^2 e^2 + 4 n b e^3 - 18 n^2 b^3 e + 6 n^2 e^2 b^2 + 4 n^3 b^3 e$$
$$- 44 n b c e a + 12 n e b c d - 8 c e^2 a - 14 e b c d + 24 n b^2 c a + 6 e^2 c d - 6 b e^3$$
$$+ 6 b^2 c d + 11 b^2 e^2 - 6 b^3 e + 12 b c e a + 3 c^2 d^2 - 6 c^2 d a - 12 c^2 n a^2 + e^4$$
$$- 6 n b^4 + 12 c^2 n^2 a^2 + 24 n^2 b c e a - 14 n b^2 c d + 12 n^3 b^2 c a)(d - 5 a + 2 a n)$$
$$(2 a n - 4 a + d) \Big/ \Big( (d + 2 a n)(d - a + 2 a n)(n-3)(6 n^2 b^2 c d + 35 n^2 b^4$$
$$+ 12 c^2 n d a + n^4 b^4 - 10 n^3 b^4 - 72 n^2 b^2 c a + 12 n c e^2 a + 70 n b^3 e - 30 n b^2 e^2$$
$$+ 4 n b e^3 - 30 n^2 b^3 e + 6 n^2 e^2 b^2 + 4 n^3 b^3 e - 92 n b c e a + 12 n e b c d$$
$$- 20 c e^2 a - 26 e b c d + 132 n b^2 c a + 6 e^2 c d - 10 b e^3 + 26 b^2 c d + 35 b^2 e^2$$
$$- 50 b^3 e + 80 b c e a + 3 c^2 d^2 - 18 c^2 d a - 36 c^2 n a^2 + e^4 - 50 n b^4 + 12 c^2 n^2 a^2$$
$$+ 24 n^2 b c e a - 26 n b^2 c d + 12 n^3 b^2 c a + 24 c^2 a^2 - 72 b^2 c a + 24 b^4) \Big) \frac{k_{n+1}}{k_n}.$$

Even more complicated results are valid in the discrete case:

$$\frac{k''_{n+1}}{k''_n} = (n+1)(2 a^2 n^3 + 12 d b n^2 - 6 a^2 n^2 + 3 d^2 n^2 + 5 a d n^2 + 12 n^2 b^2 + 24 e b n$$
$$+ 12 a e n - d^2 n - 7 a d n + 24 c a n + 12 d e n + 4 a^2 n - 12 d b n - 12 b^2 n$$
$$- 12 b e + 2 a d - 2 d^2 + 12 c d + 12 e^2)(2 a n - 2 a + d)(d - 3 a + 2 a n) \Big/ \Big($$
$$(d + 2 a n)(d - a + 2 a n)(n-1)(12 d e n + 5 a d n^2 + 12 a e n - 36 d b n$$
$$+ 12 d b n^2 + 24 e b n + 2 d^2 + 12 e^2 + 24 b d + 22 a^2 n - 12 a^2 n^2 - 12 a e$$
$$+ 3 d^2 n^2 - 7 d^2 n - 12 d e + 2 a^2 n^3 + 12 c d + 24 c a n - 36 b^2 n - 36 b e + 14 a d$$
$$+ 24 b^2 - 12 a^2 - 24 c a + 12 n^2 b^2 - 17 a d n) \Big) \frac{k_{n+1}}{k_n},$$

and

$$\frac{k'''_{n+1}}{k'''_n} = (n+1)(20\,e\,a^2\,n + 4\,b\,n^4\,a^2 - 16\,b\,n^3\,a^2 - 8\,d\,n^3\,a^2 - 2\,e\,n\,d^2 + 4\,e\,n^3\,a^2$$
$$+ 2\,d\,n^4\,a^2 - 8\,b\,a^2\,n + 4\,b\,n\,d^2 - 8\,e\,n^2\,a^2 + 12\,e\,d\,a - 13\,d^2\,n^2\,a - 4\,d\,n\,a^2$$
$$+ 5\,d^2\,n^3\,a + 10\,d\,n^2\,a^2 + 6\,d^2\,a\,n + 20\,b\,a^2\,n^2 - 36\,n^2\,b^2\,d - 8\,e\,a^2 + 12\,b\,d\,a\,n$$
$$- 8\,a\,e^2 + 16\,n\,b^3 + 8\,n^3\,b^3 + n^3\,d^3 - 24\,n^2\,b^3 - 26\,e\,a\,n\,d + 22\,e\,n^2\,a\,d$$
$$+ 10\,b\,n^3\,a\,d - 14\,b\,n^2\,d^2 - 26\,b\,n^2\,d\,a + 24\,n\,d\,b^2 - 24\,n\,b\,a\,e - 36\,n\,b\,d\,e$$
$$+ 24\,n^2\,e\,a\,b + 8\,b\,d\,e + 6\,n^2\,d^2\,e + 12\,n\,d\,e^2 + 12\,n\,d^2\,c + 6\,n^3\,d^2\,b + 24\,n\,a\,e^2$$
$$+ 24\,n\,e^2\,b + 24\,n^2\,e\,b^2 + 8\,e^3 - 8\,c\,d^2 - 24\,b\,e^2 + 16\,b^2\,e + 24\,n^2\,d\,e\,b$$
$$+ 24\,n^2\,d\,c\,a + 12\,n^3\,d\,b^2 + 24\,e\,c\,d - 48\,n\,e\,b^2 - 16\,b\,c\,d - 24\,n\,d\,c\,a + 48\,n\,e\,c\,a$$
$$+ 24\,n\,b\,c\,d - 48\,n\,b\,c\,a + 48\,n^2\,b\,c\,a - 16\,c\,e\,a - d^3\,n^2 - 4\,e\,d^2 - 2\,d^3\,n)$$
$$(d - 3\,a + 2\,a\,n)\,(2\,a\,n - 4\,a + d)\Big/\Big((d - a + 2\,a\,n)\,(d + 2\,a\,n)\,(n-2)($$
$$48\,e\,a^2\,n + 4\,b\,n^4\,a^2 - 32\,b\,n^3\,a^2 - 16\,d\,n^3\,a^2 - 14\,e\,n\,d^2 + 4\,e\,n^3\,a^2 + 2\,d\,n^4\,a^2$$
$$- 112\,b\,a^2\,n + 50\,b\,n\,d^2 - 20\,e\,n^2\,a^2 + 60\,e\,d\,a - 28\,d^2\,n^2\,a - 56\,d\,n\,a^2 + 5\,d^2\,n^3\,a$$
$$+ 46\,d\,n^2\,a^2 + 47\,d^2\,a\,n + 92\,b\,a^2\,n^2 - 24\,d^2\,a + 24\,d\,a^2 - 72\,n^2\,b^2\,d - 48\,b\,a\,d$$
$$- 40\,e\,a^2 + 94\,b\,d\,a\,n - 48\,b^3 - 32\,a\,e^2 + 88\,n\,b^3 + 8\,n^3\,b^3 + n^3\,d^3 - 48\,n^2\,b^3$$
$$+ 48\,b\,a^2 - 24\,b\,d^2 + 48\,b\,e\,a + 96\,c\,b\,a - 12\,e^2\,d + 48\,c\,d\,a - 72\,d\,b^2 - 70\,e\,a\,n\,d$$
$$+ 22\,e\,n^2\,a\,d + 10\,b\,n^3\,a\,d - 32\,b\,n^2\,d^2 - 56\,b\,n^2\,d\,a + 132\,n\,d\,b^2 - 72\,n\,b\,a\,e$$
$$- 84\,n\,b\,d\,e + 24\,n^2\,e\,a\,b + 68\,b\,d\,e + 6\,n^2\,d^2\,e + 12\,n\,d\,e^2 + 12\,n\,d^2\,c + 6\,n^3\,d^2\,b$$
$$+ 24\,n\,a\,e^2 + 24\,n\,e^2\,b + 24\,n^2\,e\,b^2 + 8\,e^3 - 20\,c\,d^2 - 48\,b\,e^2 + 88\,b^2\,e$$
$$+ 24\,n^2\,d\,e\,b + 24\,n^2\,d\,c\,a + 12\,n^3\,d\,b^2 + 24\,e\,c\,d - 96\,n\,e\,b^2 - 40\,b\,c\,d$$
$$- 72\,n\,d\,c\,a + 48\,n\,e\,c\,a + 24\,n\,b\,c\,d - 144\,n\,b\,c\,a + 48\,n^2\,b\,c\,a - 64\,c\,e\,a$$
$$- 4\,d^3\,n^2 + 4\,e\,d^2 + 3\,d^3\,n)\Big)\;\frac{k_{n+1}}{k_n}\;.$$

We omit the very lengthy result for $\frac{k''''_{n+1}}{k''''_n}$ in the discrete case. Similar expressions can be obtained for higher coefficients.
These results should indicate that the use of computer algebra is quite natural in the given context, and hand calculations seem not to be adequate for these computations.
Finally, for the sake of completeness, we present the full polynomial certificates $G_1(z,w)$ and $G_2(z,w)$ from § 4:

$$G_2(z,w) = \frac{1}{2}(-18\,w + 15\,z + 30\,z\,x^2\,w - 15\,z\,x^2\,w^2 - 10\,z^2\,w + 15\,z\,w^2 - 15\,z\,x^2 + 4\,z\,w^2\,n^3$$
$$- 15\,z\,w^2\,k + 31\,z\,w^2\,n - 16\,z\,n\,k\,x^2 + 16\,z\,n\,k - 4\,z\,n^3\,x^2 + 4\,z\,n^2\,k - 15\,z\,k\,x^2$$
$$- 31\,z\,n\,x^2 - 20\,z\,n^2\,x^2 + 15\,z\,k + 31\,z\,n + 4\,z\,n^3 + 20\,z\,n^2 - 2\,z\,w - 16\,z\,w^2\,n\,k$$
$$+ 16\,z\,w^2\,n\,k\,x^2 + 4\,z\,w^2\,n^2\,k\,x^2 - 4\,z\,w^2\,n^3\,x^2 - 4\,z\,w^2\,n^2\,k + 20\,z\,w^2\,n^2$$
$$+ 15\,z\,w^2\,k\,x^2 - 31\,z\,w^2\,n\,x^2 - 20\,z\,w^2\,n^2\,x^2 + 62\,z\,w\,n\,x^2 - 18\,z^2\,w\,n^2 - 24\,z^2\,w\,n$$
$$- 22\,w\,n^2 - 4\,w\,n^3 - 36\,w\,n + 40\,z\,w\,n^2\,x^2 + 8\,z\,w\,n^3\,x^2 - 2\,z\,w\,n + 4\,k^2\,z^2\,w\,n$$
$$- 6\,k^2\,w - 4\,z^2\,w\,n^3 + 10\,k^2\,z^2\,w - 4\,k^2\,z\,w - 4\,k^2\,w\,n - 4\,z\,n^2\,k\,x^2)\Big/k$$

and

$$G_1(z,w) = -\frac{1}{2}(5\,z^2 + 3\,z + 3\,z\,x^2\,w^2 - 3\,z\,w^2 - 3\,z\,x^2 - 32\,k\,z\,w\,n + 6\,k\,z^2\,w - 4\,z^2\,k^2\,n\,x^2$$
$$- 8\,k\,z\,w\,n^2 - 30\,k\,z\,w + 24\,k\,z\,w\,n^2\,x^2 + 78\,k\,z\,w\,x^2 - 6\,z\,w^2\,k^2 + 9\,z\,w^2\,k$$
$$- 5\,z\,w^2\,n - 10\,z^2\,k^2\,x^2 - 8\,k\,w\,n^2 - 28\,k\,w\,n + 10\,z^2\,k^2 + 15\,z^2\,k + 7\,z^2\,n - 5\,z^2\,x^2$$
$$+ 2\,z^2\,n^2 + 16\,z^2\,n\,k - 16\,z^2\,n\,k\,x^2 - 4\,z^2\,n^2\,k\,x^2 + 4\,z^2\,n^2\,k + 4\,z^2\,k^2\,n - 15\,z^2\,k\,x^2$$
$$- 7\,z^2\,n\,x^2 - 2\,z^2\,n^2\,x^2 - 4\,z\,k^2\,n\,x^2 - 12\,z\,n\,k\,x^2 + 12\,z\,n\,k + 4\,z\,n^2\,k + 4\,z\,k^2\,n$$
$$- 9\,z\,k\,x^2 - 5\,z\,n\,x^2 - 2\,z\,n^2\,x^2 - 6\,z\,k^2\,x^2 + 6\,z\,k^2 + 9\,z\,k + 5\,z\,n + 2\,z\,n^2$$
$$+ 10\,z^2\,w^2\,k^2\,x^2 + 4\,z^2\,w^2\,k^2\,n\,x^2 + 16\,z^2\,w^2\,n\,k - 16\,z^2\,w^2\,n\,k\,x^2 - 4\,z^2\,w^2\,n^2\,k\,x^2$$
$$- 10\,z^2\,w^2\,k^2 + 5\,z^2\,w^2\,x^2 - 2\,z^2\,w^2\,n^2 + 15\,z^2\,w^2\,k - 7\,z^2\,w^2\,n + 4\,z^2\,w^2\,n^2\,k$$
$$- 4\,z^2\,w^2\,k^2\,n - 15\,z^2\,w^2\,k\,x^2 + 7\,z^2\,w^2\,n\,x^2 + 2\,z^2\,w^2\,n^2\,x^2 + 6\,z\,w^2\,k^2\,x^2$$
$$+ 4\,z\,w^2\,k^2\,n\,x^2 - 5\,z^2\,w^2 + 12\,z\,w^2\,n\,k - 12\,z\,w^2\,n\,k\,x^2 - 4\,z\,w^2\,n^2\,k\,x^2 + 4\,z\,w^2\,n^2\,k$$
$$- 4\,z\,w^2\,k^2\,n - 2\,z\,w^2\,n^2 - 9\,z\,w^2\,k\,x^2 + 5\,z\,w^2\,n\,x^2 + 2\,z\,w^2\,n^2\,x^2 - 32\,k\,z^2\,w\,n\,x^2$$
$$- 8\,k\,z^2\,w\,n^2\,x^2 - 30\,k\,z^2\,w\,x^2 - 24\,k\,w + 4\,k\,z^2\,w\,n + 88\,k\,z\,w\,n\,x^2 - 4\,z\,n^2\,k\,x^2) \Big/ k$$

## Acknowlegments

The first named author thanks Tom Koornwinder and René Swarttouw for helpful discussions on their implementation `rec2ortho` [12] on the "inverse problem". Examples 2–4 given by recurrence equation (61) were provided by them. Thanks to the support of their institutions I had a very pleasant and interesting visit at the Amsterdam universities in August 1996.